\documentclass[11pt]{amsart}

\DeclareMathOperator{\im}{im} 

\DeclareMathOperator{\Div}{div}
\DeclareMathOperator{\ord}{ord}

\DeclareMathOperator{\tor}{\textup{tors}}

\newcommand{\tc}[1]{\textcolor{#1}}
\renewcommand{\phi}{\varphi}

\newcommand{\C}{\mathbb{C}}
\newcommand{\F}{\mathbb{F}}

\newcommand{\Q}{\mathbb{Q}}

\newcommand{\Z}{\mathbb{Z}}

\usepackage[utf8]{inputenc}
\usepackage[T1]{fontenc}
\usepackage{amsmath}
\usepackage{amsthm}
\usepackage{mathtools}  
\usepackage[dvipsnames]{xcolor}     
\usepackage{graphicx}               
\usepackage{verbatim}               
\usepackage{enumerate,enumitem}     
\usepackage[charter]{mathdesign}              
\usepackage{geometry}   
    \newlength{\alphabet}
\usepackage[]{todonotes}          
             
\usepackage{longtable}             
\usepackage{lscape}              
\usepackage{makecell}           
\usepackage{breqn}                  

\usepackage[colorlinks=true, urlcolor=RoyalBlue, citecolor=Green, linkcolor=purple]{hyperref}   

\graphicspath{{img/}{../img/}}                
\DeclareGraphicsExtensions{.pdf, .png, .JPG}    

\theoremstyle{definition}
\newtheorem{definition}{Definition}[section]
\newtheorem{notation}[definition]{Notation}
\newtheorem{example}[definition]{Example}
\newtheorem{remark}[definition]{Remark}

\theoremstyle{plain}

\newtheorem*{conjecture*}{Conjecture}
\newtheorem{conjx}{Conjecture}
 
\newtheorem*{theorem*}{Theorem}

\newtheorem*{corollary*}{Corollary}

\newtheorem{theorem}[definition]{Theorem}
\newtheorem{proposition}[definition]{Proposition}
\newtheorem*{proposition*}{Proposition}
\newtheorem{claim}[definition]{Claim}
\newtheorem*{claim*}{Claim}
\newtheorem{lemma}[definition]{Lemma}
\newtheorem*{note*}{\tc{teal}{Note}}
\newtheorem{corollary}[definition]{Corollary}

\newcommand{\Gal}{\operatorname{Gal}}
\newcommand{\SL}{\operatorname{SL}}
\newcommand{\GL}{\operatorname{GL}}

\newcommand{\cusp}{\operatorname{Cusp}}

\newcommand{\CC}{\mathcal{C}}

\newcommand{\p}{\mathfrak{p}}
\newcommand{\n}{\mathfrak{n}}
\newcommand{\m}{\mathfrak{m}}
\newcommand{\fa}{\mathfrak{a}}

\newcommand{\fd}{\mathfrak{d}}

\newcommand{\pp}{|\p|}

\usepackage[mathcal]{euscript}
\usepackage{enumitem}
\usepackage[dvipsnames]{xcolor}
\usepackage{tikz-cd}
\usepackage{geometry}
\title{Rational torsion of generalised Drinfeld modular Jacobians of prime power level}

\author{Mar Curc{\'o} Iranzo}
\thanks{I am grateful to Valentijn Karemaker and Gunther Cornelissen for discussions on this topic and comments
on an earlier version of the paper. I would also like to thank the reading committee of my PhD thesis for their comments on the version of this manuscript written in there.}
\email{m.curcoiranzo@uu.nl}
\address{Mathematical Institute, Utrecht University, PO Box 80010, 3508 TA, Utrecht, the Netherlands}
\date{\today}

\begin{document}

\begin{abstract}
For a prime $\p \subseteq \F_{q}[T]$ and a positive integer $r$, we consider the generalised Jacobian $J_{0}(\n)_{\mathbf{m}}$ of the Drinfeld modular curve $X_{0}(\n)$ of level $\n=\p^r$, with respect to the modulus~$\mathbf{m}$ consisting of all cusps on the modular curve. We show that the $\ell$-primary part of the group $J_{0}(\n)_{\mathbf{m}}(\F_{q}(T))_{\rm{tor}}[\ell^{\infty}]$ is trivial for all primes $\ell$ not dividing $q(q^{2}-1)$. Our results stablish an analogue to those of Yamazaki--Yang \cite{yamazaki2016rational} for the classical case. 
\end{abstract}

\maketitle

\section{Introduction}

\subsection{Torsion of modular (semi-)abelian varieties.}
Given a positive integer $N$, let $X_{0}(N)$ denote the modular curve of level $N$ for the congruence subgroup $\Gamma_{0}(N)=\{\begin{psmallmatrix} a & b \\ c & d \end{psmallmatrix}\in \rm{SL}_{2}(\Z) : c \equiv 0 \pmod{N}\}$. This is a projective nonsingular curve defined over $\Q$. 
For $N$ prime, Ogg \cite[Conjecture~2]{ogg1975diophantine} considered the Jacobian variety $J_{0}(N)$ of $X_{0}(N)$ as a natural family of candidates of higher-dimensional abelian varieties whose torsion could be determined. His predictions on $J_{0}(N)$, based on computational evidence, are now known as Ogg's Conjecture. 
To state it, write $C_{N}$ for the image of the degree-zero divisors with support at the cusps of $X_{0}(N)$ in $J_{0}(N)$ and define $C_{N}(\Q) \coloneqq C_{N} \cap J_{0}(N)(\Q)$; Ogg proved that, 
for $N\geq 5$ prime, the cyclic group $C_{N}(\Q)$ has order $n=\frac{p-1}{\gcd(p-1, 12)}$ and made the following conjecture: the cyclic  group $C_{N}(\Q)$ of  order  $n$ is  the full torsion  subgroup $J_{0}(N)(\Q)_{\tor}$ of $J_0(N)(\Q)$.
The conjecture was proved later by Mazur \cite[Theorem~I]{mazur1977modular} using the Eisenstein ideal machinery developed in his work. 
A natural generalisation of this conjecture is the following. 
\begin{conjx}[Generalised Ogg's conjecture
] \label{conj:A} Let $N$ be a positive integer. Then
\begin{equation*}
    C_{N}(\Q) = J_0(N)(\Q)_{\tor}.
\end{equation*}
\end{conjx}
Notice that by theorems of Manin \cite[Corollary~3.6]{manin1972parabolic} and Drinfeld \cite{drinfeld1973two} the group $C_{N}$ is finite, thus, the inclusion  $C_{N}(\Q) \subseteq J_0(N)(\Q)_{\tor}$ always holds. Moreover, if $C(N)$ is the subgroup formed by the equivalence classes of $\Q$-rational divisors in $C_{N}(\Q)$, i.e., the divisor classes fixed under the action of the absolute Galois group of~$\Q$, a further conjecture states that $C_{N}(\Q)$ is equal to $C(N)$. 

\begin{conjx}[Rational cuspidal groups conjecture] \label{conj:B} Let $N$ be a positive integer. Then
\begin{equation*} 
    C_{N}(\Q) = C(N).
\end{equation*}
\end{conjx}

A strategy for proving Conjectures \ref{conj:A} and \ref{conj:B} consists in decomposing the groups $C(N), C_{N}(Q)$ and $J_{0}(\Q)_{\tor}$ into the direct sum of their $\ell$-primary parts for all primes $\ell$ and showing that each respective $\ell$-primary part equals that of the other groups. In this spirit, the strongest result to date says that $C(N)[\ell^{\infty}]$ and $J_{0}(\Q)_{\tor}[\ell^{\infty}]$ agree for all primes $\ell$ with $\ell \geq 5$, or $\ell=3$ and either $3^{2} \nmid N$ or $N$ has a prime divisor congruent to $-1$ modulo $3$. This result is due to Yoo \cite{yoo2021rational} and is based on previous work of Ohta \cite{ohta2013eisenstein} on the Eistenstein ideal. See Section 1 in \cite{yoo2021rational} for a summary of more (partial) results on Conjectures~\ref{conj:A} and \ref{conj:B}.
On the other hand, the elaborate task of determining the structure of $C(N)$ has recently been accomplished in further work by Yoo~\cite{yoo2019rationalcusp}. A meticulous study of modular units, in particular of $\eta$-quotients, allows him to determine the structure and generators of $C(N)[\ell^{\infty}]$ for all $\ell$. Hence, as of today we can say that the structure of 
 $J_{0}(\Q)_{\tor}[\ell^{\infty}]$ is understood for the primes $\ell$ in Yoo's result \cite[Theorem 1.4]{yoo2021rational}.

In a parallel direction, a number of authors have moved on to determining the structure of $J_{0}(N)_{\mathbf{m}}(\Q)_{\tor}$, where $J_{0}(N)_{\mathbf{m}}(\Q)_{\tor}$ is the generalised Jacobian of $X_{0}(N)$ with respect to a cuspidal modulus $\mathbf{m}$. That is, if we let $\cusp(X_{0}(N))$ be the set of cusps of $X_{0}(N)$ and $\mathbf{m}= \sum_{\omega \in \cusp(X_{0}(N))} \omega$, the generalised Jacobian is the group of invertible line bundles on the singular curve that results from $X_{0}(N)$ after the identification of all the points in $\mathbf{m}$. The general construction of this group variety can be found in \cite{rosenlicht1954generalized} and \cite{serre1959groupes}. In our case, $J_{0}(N)_{\mathbf{m}}(\Q)_{\tor}$ is what is called a semi-abelian variety; that is, the extension of $J_{0}(N)$ -- an abelian variety -- by a commutative affine linear algebraic group $L_{\mathbf{m}}$ isomorphic to a number of copies of $\mathbb{G}_{m}$. 
The group $J_{0}(N)_{\mathbf{m}}(\Q)$ is not not necessarily finitely generated anymore, but the subgroup $J_{0}(N)_{\mathbf{m}}(\Q)_{\tor}$ is finite. Just like the Jacobian $J_{0}(N)$, the generalised Jacobian $J_{0}(N)_{\mathbf{m}}$ is an object of high (arithmetic) relevance when it come to the study of the curve $X_{0}(N)$. With this underlying motivation, Yamazaki-Yang \cite{yamazaki2016rational}, Wei-Yamazaki~\cite{wei2019rational}, Curc\'o-Iranzo~\cite{iranzo} and Jordan-Ribet-Scholl~\cite{JordanRibetScholl} study the structure of $J_{0}(N)_{\mathbf{m}}(\Q)_{\tor}$ for different $N$ and with different approaches. The following result of Yamazaki-Yang is relevant for us in this paper. 

\begin{theorem}\cite[Theorem 1.1.3]{yamazaki2016rational} \label{thm:yamazakiyang}
Let $p$ be a prime number and $n$ be a positive integer. For $\ell \nmid 2p$, the group $J_{0}(p^{n})_{\mathbf{m}}(\Q)[\ell^{\infty}]$ is trivial. Moreover, if Conjectures~\ref{conj:A} and \ref{conj:B} hold, we have 
$$ J_0(p^n)_{\mathbf{m}}(\Q)_{\tor} \simeq\prod_{i=0}^{\lfloor\frac{n-1}{2}\rfloor}
\Z/(2p^{i}\Z) \times \prod_{i=1}^{\lfloor \frac{n}{2}\rfloor}
\Z/(2p^{i}\Z).$$
\end{theorem}

\subsection{Torsion of Drinfeld modular (semi-)abelian varieties.} In the current paper, we study an analogue of the results of Yamazaki-Yang in the function field setting. 

\begin{notation} Let $p$ be a prime and fix $q=p^{t}$. We fix the following notation for the rest of the paper. 

\begin{tabular}{p{0.08\textwidth}p{0.8\textwidth}}
$\mathbb{F}_q$ & = finite field of characteristic $p$ with $q$ elements\\
$\mathcal{A}$ & = $\mathbb{F}_q[T]$ polynomial ring in $T$ over $\mathbb{F}_q$\\
$K$ & = $\mathbb{F}_q(T)$ rational function field\\
$K_\infty$ & = $\mathbb{F}_q((\pi))$ the completion of $K$ at the infinite place; $\pi:=T^{-1}$\\
$|\cdot|$ & = $|\cdot|_\infty$ = normalised absolute value on $K_\infty$; that is, $|T|_\infty:=q$\\
$\mathcal{O}_\infty$ & = $\mathbb{F}_q[[\pi]]$ ring of integers in $K_\infty$\\
$\mathbb{C}_\infty$ & = the completion of an algebraic closure of $K_\infty$\\
$G$ & = group scheme $\GL(2)$ over $\mathbb{F}_q$\\
$Z$ & = scalar matrices in $G$\\
\end{tabular}
\end{notation}

Before introducing our main result, we present the analogue of Conjectures~\ref{conj:A} and \ref{conj:B} for Drinfeld modular curves. Given an monic element $\n\in \mathcal{A}$, let $X_{0}(\n)$ denote the Drinfeld modular curve of level $\n$ for the congruence subgroup $\Gamma_0(\n)$ consisting of matrices in $\Gamma = G(\mathcal{A})$ that are upper triangular modulo $\n$ -- see Section~\ref{sec:background}. We regard $X_{0}(\n)$ as a smooth projective curve defined over $K$. Let $J_0(\n)$ be the Jacobian variety of $X_0(\n)$ and $J_0(\n)(K)_{\text{tors}}$ its rational torsion subgroup. The Lang-N\'eron theorem \cite{LangNeron}, which generalises the Mordell-Weil theorem, asserts that 
$J_0(\n)(K)$ is finitely generated, and thus that $J_0(\n)(K)_{\text{tors}}$ is finite. Similar as before, let $\CC_\n$ be the cuspidal subgroup of $J_{0}(\n)$, that is, the image of the degree-zero divisors with support at the cusps of $X_{0}(\n)$ in $J_{0}(\n)$, and let $\CC_\n(K) := \CC_\n \cap J_{0}(\n)(K)$ be the rational cuspidal subgroup of $J_0(\n)$. Moreover, let
$\CC(\n)$ be the rational cuspidal divisor class group of $X_0(\n)$, which is the subgroup of $\CC_{\n}$ formed by the equivalence classes of $K$-rational divisors in $C_{N}(\Q)$, i.e., the divisor classes fixed under the action of the absolute Galois group of~$K$. Since $\CC_{\n}$ is a finite group -- it is generated by a finite set of torsion elements \cite{GekelerFiniteness} --, we have the tower of inclusions 
\begin{equation}\label{eq:incl}
    \CC(\n)\subseteq \CC_\n(K)\subseteq J_{0}(\n)(K)_{\tor}.
\end{equation}
Analogues of Conjectures~\ref{conj:A} and \ref{conj:B} in this setting predict the equality of two the inclusions. 

\begin{conjx} \label{conj:C} Let $\n \in \mathcal{A}$ be monic. Then we have
\begin{equation*}
    \CC(\n)= \CC_\n(K)=J_{0}(\n)(K)_{\tor}.
\end{equation*}
\end{conjx}

As analogous results to those of Ogg and Mazur in this direction, Gekeler~\cite{GekelerDiscriminant}
computed the order of the cyclic group $\CC_{\p}$ and P\'al \cite{Pal} proved Conjecture \ref{conj:C} for $\n=\p$ prime. Their techniques resembled those of the classical modular setting, that is, they employed a study of the discriminant function and the development of the appropriate Eisenstein ideal machinery in this setting. Further progress in determining the structure of $\CC(\n)$ and proving Conjecture \ref{conj:C} in full generality can be found in works of Papikian-Wei~\cite{PapikianWei} and Ho \cite{kevincusp},  \cite{HoTorsion}. See \cite{armanaoggfunctionfields} for a survey on the topic.

In the Drinfeld case, one can also consider the generalised Jacobian $J_{0}(\n)_{
\mathbf{m}}$ of $X_{0}(\n)$ with cuspidal modulus $\mathbf{m} = \sum_{\omega \in \cusp(X_{0}(\n))} \omega$. In \cite{wei2019rational}, Wei-Yamazaki compute the structure of  $J_{0}(\n)_{
\mathbf{m}}(K)_{\tor}$ up to $q(q^{2}-1)$-torsion. The goal of this paper is two-fold: compute the structure of the group $J_{0}(\n)$ for $\n=\p^{r}$ and  compare the results for this level to the ones in the classical case when $N=p^{r}$ given in \ref{thm:yamazakiyang}. Our main result is the following. 

\begin{theorem} \label{thm:mainthm}
Let $\p \in \mathcal{A}$ monic irreducible and $r$ an integer. For all odd primes $\ell$ with $\ell \nmid q(q^{2}-1)$ the subgroup $J_{0}(\p^{r})_{\mathbf{m}}(K)[\ell^{\infty}]$ is trivial. Moreover, if Conjecture \ref{conj:C} holds, we have $J_{0}(\p^{r})_{\mathbf{m}}(K)_{\tor} \simeq \prod_{i=0}^{r-1} 
\Z/(q-1)\Z$. 

\end{theorem}

The results in Theorem~\ref{thm:mainthm} show that the situation in the Drinfeld case is analogous to the one in Theorem~\ref{thm:yamazakiyang}. 

The structure of the paper is the following. In Section~\ref{sec:background} we present some background on the Drinfeld modular curves $X_{0}(\n)$ and their cusps. In particular, in Theorem~\ref{thm:orbits} we describe the residue fields of the cuspidal closed points of $X_{0}(\n)$. We also review known results on $J_{0}(\n)$ and $\CC(\n)$ in the  direction of Conjecture~\ref{conj:C}. In Section~\ref{sec:generalisedJacobians} we recall some results on the Drinfeld generalised Jacobian $J_{0}(\n)_{\mathbf{m}}$. We introduce a group homomorphism denoted 
by $\delta$ and reduce the proof of Theorem~\ref{thm:mainthm} to determining the injectivity of a related homomorphism $\overline{\delta}$ on $\CC(\n)$. Finally, in the Sections~\ref{sec:imagedelta} and \ref{sec:kerneldelta} we prove the main result. The strategy of the proof then goes as follows. In Section~\ref{sec:imagedelta} we develop the tools to compute the image of certain cuspidal divisors appearing in Theorem~\ref{thm:cuspgen}. In Section~\ref{sec:kerneldelta} we use that computation to construct the matrix that represents the map $\overline{\delta}|_{\CC(\p^r)}$. Finally, 
that representation allows us to see that $\overline{\delta}|_{\CC(\p^r)}$ is injective.

\section{Background} \label{sec:background}

Let $\Omega = \mathbb{C}_\infty - K_\infty$ be the Drinfeld upper half plane. 
The group $\Gamma_0(\n)$ acts on $\Omega$ by linear fractional transformations. The Drinfeld modular curve $Y_0(\n)$ can be described analytically as the quotient $ \Gamma_0(\n)\backslash \Omega$. The points of $Y_{0}(\n)$ parametrise rank 2
Drinfeld modules equipped with a cyclic subgroup of order $\n$.
Then, $X_{0}(\n)$ is the unique smooth projective curve over~$K$ containing $Y_0(\n)$, that is, the canonical non-singular compactification of $Y_{0}(\n)$. The set of cusps of $X_{0}(\n)$ is $\cusp(X_{0}(\n)) = X_{0}(\n) \setminus Y_{0}(\n)$. after base change to $\C_{\infty}$, there is a canonical 1-to-1 correspondence between cusps of $X_{0}(\n)_{\C_{\infty}}$ and 
$\Gamma_{0}(\n) \backslash \mathbb{P}_{1}(K)$. 
A divisor $D$ of $X_{0}(\n)_{\C_{\infty}}$ with  $\text{supp}(D) \subseteq \text{Cusps}(X_{0}(\n)_{\C_{\infty}})$ is called a cuspidal divisor.
Then the group $\CC_\n$ (resp. $\CC(\n)$) defined in the introduction is the subgroup of $J_{0}(\n)$ of classes of cuspidal divisors (resp. $K$-rational cuspidal divisors).
A recent remarkable result on Conjecture~\ref{conj:C} is the following. 

\begin{theorem}\cite[Theorem 1.7] {HoTorsion}\label{thm:kevintorsion} Let  $\n = \p^{r} \in \mathcal{A}$ be a prime power ideal with $r \geq 2$. For any prime $\ell \nmid  q(q - 1)$, we have
$$\CC(\p^r)[\ell^{\infty}] = J_{0}(\p^{r})(K)_{\tor}[\ell^{\infty}].$$
\end{theorem}

Motivated by Conjecture~\ref{conj:C}, a number of authors have studied the structure of the group $\CC(\n)$ in order to understand the torsion of $J_{0}(\n)$ better. For that, understanding the points in $\text{Cusps}(X_{0}(\n))$ is crucial. We proceed by giving some background knowledge about cusps. 

\begin{theorem}\textup{(\cite[pg. 196]{GekelerDiscriminant}, \cite[Lemma 3.1(1)]{PapikianWei})} \label{thm:orbits}
Let $\n \in \mathcal{A}$ be monic. The cusps of $X_0(\n)_{\C_{\infty}}$ are in
bijection with $\Gamma_0(\n) / \mathbb{P}^1(K)$. Every cusp in $\cusp(X_{0}(\n))_{\C_{\infty}}$ has a representative $\left[\begin{smallmatrix} \fa \\ \m \end{smallmatrix}\right]$ for some elements $\fa, \m \in \mathcal{A}$, with $\m \mid \n$ and $\gcd(\fa , \n)=1$.
In particular, we can describe them as the set
\[ \cusp(X_{0}(\n)_{\C_{\infty}}) \simeq \left\{\left[\begin{smallmatrix} \fa \\ \m\end{smallmatrix}\right] \in \Gamma_0(\n)\setminus\mathbb{P}^1(K): \fa, \m \in \mathcal{A} \text{ monic, } \m\mid \n, \text{ and } \gcd(\fa, \n) = 1 \right\}.\]
For two divisors $\m$ and $\m'$ of $\n$, we have
    \[
    \left[\begin{smallmatrix} \fa \\ \m \end{smallmatrix}\right]=\left[\begin{smallmatrix} \fa' \\ \m' \end{smallmatrix}\right] \iff \m=\m' \text{ and } k \fa \equiv \fa' \pmod{\fd},
    \]
    where $\fd= \fd(\m):=\gcd(\m, \n/\m)$ and $k \in \F_{q}^{\times}$.
\end{theorem}

A cusp $\omega$ of $X_0(\n)_{\C_{\infty}}$ with a representative of the form $\left[\begin{smallmatrix} \fa \\ \m \end{smallmatrix}\right]$ as in Theorem~\ref{thm:orbits} is said to be of height $\m$ \cite[(2.4)]{gekelerfundamental}. Let $(P_{\m})$ be the divisor given by the sum of all the cusps of $X_0(\n)_{\C_{\infty}}$ of fixed height $\m \mid \n$. The cuspidal divisors $(P_{\m})$ are $K$-rational, that is, they are invariant under the action $\text{Gal}(\overline{K}/K)$, cf. \cite[Prop 6.3]{GekelerInvariants}. Indeed, the cusps of
$X_0(\n)_{\C_{\infty}}$ of the same height form an orbit under $\text{Gal}(\overline{K}/K)$.  For $\m=\p^{i}$ with $i \in \{0, \dots, r\}$, we let $P_{i}$ be the closed $K$-point of $X_{0}(\n)/K$ given by $P_{i}:=(P_{\m})$. The following results describe the residue fields of the closed points $P_{i}$.

\begin{proposition} \label{prop:orbitsize} \label{prop:orbitSL2} Let $\n = \p^{r}$ and let $\m \mid \n$. Then the orbit of the cusp $\left[ \begin{smallmatrix} 1 \\ \m \end{smallmatrix}\right]$ under the action of $\textup{GL}_{2}(\mathcal{A})$ is of size $\varphi(\m, \n/\m)/(q-1)$. 
\end{proposition}

\begin{proof}
    This follows from Theorem~\ref{thm:orbits}.
\end{proof}

\begin{definition}
    For a polynomial $M\in \F_{p}[T]$, let $K(M)$ denote the cyclotomic function field obtained after adjoining the $M$-torsion of the Carlitz module in $\overline{\F_{p}(T)}$ to $K=\F_{p}(T)$. See \cite{conradcarlitz} for a survey on cyclotomic function fields.
\end{definition}

\begin{theorem} \label{thm:field} Let $\m=\p^{i}$, $\fd=\fd(\m)=\gcd(\m, \n/\m)=\p^{\min\{i, r-i\}}$, and let $P_{i}$ be the closed point defined by the cuspidal divisor $(P_{\m})$. The residue field $K(P_{i})$ of the point $P_{i}$ is $K^{+}(\mathfrak{d})$ the maximal real subfield of the $\mathfrak{d}$-th cyclotomic field $K(\mathfrak{d})$ of $K$; that is, the fixed field by the automorphisms of $\Gal(K(\mathfrak{d})/K)$ induced by elements of $\F_{q}^{\ast}$. See \cite{conradcarlitz} for a survey on cyclotomic function fields.
\end{theorem}

\begin{proof}
    By classical results of Gekeler~\cite[Rmk 4.7]{GekelerInvariants}, we know that the cusps of $X_{0}(\n)(\C_{\infty})$ are $K^{+}(\n)$-rational points. Hence, we know that $K(P_{i})$ is a subfield of the cyclic field $K^{+}(\n)$. 
    Consider now the cusp $\left[ \begin{smallmatrix} 1 \\ \m \end{smallmatrix}\right]$. By Proposition~\ref{prop:orbitSL2}, we know that the orbit of $\left[ \begin{smallmatrix} 1 \\ \m \end{smallmatrix}\right]$ under the action of $\text{GL}_{2}(\mathcal{A})$ is of size $\varphi(\m, \n/\m)/(q-1)$. 
    It follows that the field $K(P_{i})$ is the unique degree $\varphi(\m, \n/\m)/(q-1)$ subfield of the cyclic extension $K^{+}(\n)/K$. That is, we get that $K(P_{i})=K(\mathfrak{d})^{+}$. \end{proof}

In particular, there are two points in $\cusp(X_{0}(\n))$ defined over $K$: the cusp at 0 (given by the class $\left[  \begin{smallmatrix}
    0 \\ 1
\end{smallmatrix}\right]$) and the cusp at $\infty$ (given by the class $\left[ \begin{smallmatrix}
    1 \\ \n
\end{smallmatrix}\right]$). We denote them by $\omega_{0}$ and $\omega_{\infty}$ respectively. 

When $\n=\p$ is a monic irreducible polynomial $\p \in \mathcal{A}$, $\omega_{0}$ and $\omega_{\infty}$ are the only cusp points of $X_{0}(\n)$. Hence, the group $\CC(\n)$ is the cyclic group generated by the class of the divisor $\omega_{0}- \omega_{\infty}$. Gekeler shows in  \cite{GekelerDiscriminant} that the latter is an element of order $N(\p)$, see Definition~\ref{def:order}. When the level $\n$ is not prime, determining the structure of $\CC(\n)$ is more involved. In \cite{kevincusp}, Ho studies the structure of the group $\CC(\n)$ for $\n=\p^{r}$ a power of a monic irreducible polynomial $\p \in \mathcal{A}$. In particular, he constructs a basis of $\CC(\p^{r})$ for any $r \geq 1$. The constructed basis will be very useful to prove our results. For convenience of the reader, we reproduce the results of Ho below.  

\begin{definition} \label{def:order} Fix a prime $\p \in \mathcal{A}$. We define the integers $M(\p)$ and $N(\p)$ as follows: 
\[ M(\p) \coloneqq \frac{\pp^{2}-1}{q^{2}-1} \quad \text{and} \quad N(\p):=\begin{cases}
\frac{\pp-1}{q^2-1},  & \text{if $\deg(\p)$ is even;} \\
\frac{\pp-1}{q-1},  & \text{otherwise.}
\end{cases}\]
    
\end{definition}

\begin{definition} \cite[Section 3.2]{kevincusp} \label{defn:definitiondivisors} For each $i \in \{1, \ldots, r-2\}$ we define the divisor $C_{i} \coloneqq P_{i} - \deg(P_{i}) [\infty] \in \text{Div}_{\text{cusp}}^{0}(X_{0}(\n)(K))$. We further define $D_{0}$ and $D_{r-1} \in \text{Div}_{\text{cusp}}^{0}(X_{0}(\n)(K))$ as 
\[D_{0}\coloneqq C_{0} +(q-1)\left(\sum_{1 \leq i \leq \lfloor\frac{r}{2}\rfloor}C_{i} + \sum_{\lfloor\frac{r}{2}\rfloor+1 \leq i \leq r-1} \pp^{2i-r}C_{i}\right)\]
 and $D_{r-1}$ through the case-by-case definition given in the Appendix -- see Definition~\ref{def:Dr-1} -- for convenience of the reader. 
\end{definition}

\begin{theorem}\cite[Theorem 3.5]{kevincusp} \label{thm:cuspgen}
    Let $\p\in A$ be a prime and $r\geq 2$. The elements $C_i$, $D_{r-1}$, and $D_0$ are defined in Definition~\ref{defn:definitiondivisors}. Then $$\CC(\p^r) = \left(\bigoplus_{i=1}^{m}\langle \overline{C_i}\rangle\right)\oplus\left(\bigoplus_{i=m+1}^{r-2}\langle \overline{C_i-\pp C_{i+1}}\rangle\right)\oplus \langle \overline{D_{r-1}}\rangle\oplus \langle  \overline{D_0}\rangle,$$ where $m:=\lfloor \frac{r-1}{2}\rfloor$ and
\begin{enumerate}
    \item $\ord(\overline{C_i})=\pp^{r-i}M(\p)$ for $1\leq i\leq m$;
    \item $\ord(\overline{C_i-\pp C_{i+1}})=\pp^{i}M(\p)$ for $m+1\leq i\leq r-2$;
    \item $\ord(\overline{D_{r-1}})=M(\p)$;
    \item $\ord(\overline{D_0})=N(\p)$.
\end{enumerate}
\end{theorem}

Theorem~\ref{thm:cuspgen} is obtained after taking a close look at the modular units of $X_{0}(\n)$ for $\n=\p^{r}$. In particular, in \cite{kevincusp} the author constructs a map $g : \text{Div}^{0}(\cusp(X_{0}(\n)) \longrightarrow  \mathcal{E}_{\n} \otimes_{\Z} \Q$
which, in analogy to the classical case, to each degree zero cuspidal divisor $D$ attaches a modular unit that has $D$ as a divisor. This is done through so-called $\Delta$-quotients. 

\begin{definition} Let $\n \in \mathcal{A}$ and consider the Drinfeld discriminant function $\Delta(z) : \Omega \rightarrow \C_{\infty}$ defined by Gekeler in \cite{GekelerDiscriminant}. For $\m \in \mathcal{A}$ monic, let $\Delta_{\m}:= \Delta(\m z)$. All functions $\Delta_{\m}$ with $\m \mid \n$ are modular forms on $\Omega$ of weight $q^{2}-1$ and type 0 for $\Gamma_{0}(\n)$, see \cite{GekelerDiscriminant}. A $\Delta$-quotient on $\Gamma_{0}(\n)$ is a function in $\mathcal{E}_{\n}\otimes_{\Z}\Q$ for $\Gamma_{0}(\n)$ of the form $\prod_{\m\mid\n} \Delta_{\m}(z)^{r_{\m}}$, with $r_{\m} \in \Q$ and $\sum_{\m \mid \n}r_{\m} =0$.  
\end{definition}

\begin{remark}
    Notice that since $\frac{\Delta}{\Delta_{\m}} \in \mathcal{E}_{\n}$ for all $\m \mid \n$, a $\Delta$-quotient is indeed an element of $\mathcal{E}_{\n}\otimes_{\Z} \Q$. 
\end{remark}

The map $g$ can then be constructed by understanding the order of the zeroes of the functions $\Delta_{\m}$ at all the cusps of $X_{0}(\n)$. In \cite{kevincusp}, the author uses the orders of such functions for $\n=\p^{r}$ described in \cite{GekelerDiscriminant} and constructs  
\[\begin{array}{cc}
    &g : \text{Div}^{0}(\cusp(X_{0}(\n))(K) \quad \longrightarrow  \mathcal{E}_{\n} \otimes_{\Z} \Q \\
    \\
     &C:= \displaystyle \sum_{\makecell{\m\mid \n \\ \text{monic}}} a_{\m}P_{\m} \quad \quad \longmapsto \quad   \prod_{\makecell{\m\mid \n \\ \text{monic}}} \Delta_{\m}(z)^{r_{\m}}
\end{array}\]
where the linear map $g$ is described by 
\begin{equation}
    \begin{pmatrix} r_{1} \\ r_{\p} \\
    \vdots \\ 
    r_{\p^{r}} \end{pmatrix} = \frac{1}{(q-1)\pp^{r-1}(\pp^2-1)} \Upsilon(\p^{r}) \cdot \begin{pmatrix} a_{1} \\ a_{\p} \\
    \vdots \\ 
    a_{\p^{r}}
    \end{pmatrix},
\end{equation}
and the matrix $\Upsilon(\p^{r}) = \left(\upsilon[i][j]\right)_{i, j}$ is described by 
$$ 
\upsilon[i][j] = \begin{cases}
(\pp^{2}+1)\pp^{m(j)-1}, & \text{if } 1\leq  i = j \leq r-1; \\  
-\pp^{m(j)},  & \text{if } |i-j|=1 \text{ and } j \neq 0, r; \\
(q-1)\pp, & \text{if } (i, j)=(0, 0) \text{ or } (r, r); \\
1-q, & \text{if } (i, j)=(1, 0) \text{ or } (r-1, r); \\
0, & \text{otherwise.}\\
\end{cases}
$$
Then, the strategy to prove Theorem~\ref{thm:cuspgen} is to use an appropriate root of the $\Delta$-quotient $\frac{\Delta(z)}{\Delta(\n z)}$ to to find an optimal way of computing the order of a cuspidal divisor $D$ in $\mathcal{C}(\n)$. 
This is also used to determine whether two $\Delta$-quotients have associated cuspidal divisors that generate groups in $J_{0}(\n)$ with trivial intersection. This allows Ho to determine a cyclic decomposition of $\CC(\n)$.

The functions $g(D)\in \mathcal{E}_{\n}\otimes_{\Z}\Q$ for $D \in \{C_{i}: 1 \leq i \leq \lfloor{\frac{r-1}{2}\rfloor}\}\cup\{C_{i}-\pp C_{i+1}:\lfloor{\frac{r+1}{2}\rfloor} \leq i \leq r-2\}\cup\{D_{r-1}\}$ will play an important role in proving our results. In particular, we will need a way of computing their leading Fourier coefficients at the different cusps of $X_{0}(\n)$. The following result provides such a tool and will be used later in Proposition~\ref{prop:coef}.

\begin{lemma} \label{lem:identityDelta} Let $\left[\begin{smallmatrix}
    \fa' \\ \m'
\end{smallmatrix}\right] \in \cusp(X_{0}(\n)(\C_{\infty}))$ be a cusp  with $\m'=\p^{k}$ and $\gcd(\fa', \m')=1$. Let $A =\begin{pmatrix} a & b \\ \p^{k} & \fa' \end{pmatrix} \in \mathrm{SL}_{2}(\mathcal{A})$ 
be a matrix such that $A^{-1} \cdot \omega_{\infty} = \left[\begin{smallmatrix}
    \fa' \\ \m'
\end{smallmatrix}\right]$. Then, for a monic divisor $\m \mid \n$ we have the identity
$$ \Delta_{\m}(A^{-1}\cdot z)= \left(\frac{d}{\m}(-\p^{k}z+a)\right)^{q^{2}-1} \Delta\left(\frac{d^{2}}{\m}z+\frac{\eta d}{\m}\right),$$
where $d = \gcd(\m, \p^{k})$ and $\eta \in \mathcal{A}$. 
\end{lemma}

\begin{proof}
By definition we have 
\[ \Delta_{\m}(A^{-1}z)= \Delta\left(\frac{\m\fa'z-\m b}{\p^{k}z+a}\right) = \Delta(\alpha z),  \]
where $\alpha = \begin{psmallmatrix} \m\fa' & -\m b \\ \p^{k} & a \end{psmallmatrix}.$ 
Take $L = \begin{psmallmatrix}
     x & y \\ u & v
\end{psmallmatrix} \in \SL_{2}(\mathcal{A})$ 
with $x= \frac{\m\fa'}{\gcd(\p^{k}, \m\fa')}=\frac{\m \fa'}{d}$ and $u=\frac{-\p^{k}}{d}$. 
Then we have 
$L^{-1}\alpha = \begin{psmallmatrix} \ast & \ast \\ 0 & \m/d
\end{psmallmatrix}$. And since $\text{det}(L^{-1}\alpha) = \text{det}(L^{-1}) \text{det}(\alpha) = 1 \cdot \m = \m$, we infer that $L^{-1}\alpha = \begin{psmallmatrix} d & \eta \\ 0 & \m/d
\end{psmallmatrix}$, where $\eta$ is some element in $\mathcal{A}$. Since $\Delta(z)$ is a modular form of weight $q^2-1$ under elements of $\SL_{2}(\mathcal{A})$, we get the following:
$$ \Delta_{\m}(A^{-1}z)= \Delta(\alpha z) = \Delta(LL^{-1}\alpha z) $$ 
$$ =\left(u\left(\frac{dz+\eta}{\m/d} + v\right)z+v\right)^{q^{2}-1} \Delta\left(\frac{dz+\eta}{\m/d}\right)=\left(\frac{d}{\m}(-\p^{k}z+a)\right)^{q^{2}-1} \Delta\left(\frac{d^{2}}{\m}z+\frac{\eta d}{\m}\right).$$ 

\end{proof}

\section{The generalised Jacobian} \label{sec:generalisedJacobians}
Let $J_{0}(\n)_{\mathbf{m}}$ be the generalised Jacobian of $X_{0}(\n)$ with respect to the modulus $\mathbf{m}:= \sum_{P \in C} P$, where $C:=\text{Cusps}(X_{0}(\n))$. That is, we have
\[
J_{0}(\n)_{\mathbf{m}}(K) = \text{Div}^{0}(X_{0}(\n) \setminus C)/\{\text{div}(f) \mid  f \in K(X_{0}(\n))^{\times}, f\equiv  1 \pmod{\mathbf{m}}\}.
\]
Hence, $J_{0}(\n)_{\mathbf{m}}$ fits into the exact sequence
\begin{equation} \label{eq:generalisedJacobian}
    0 \rightarrow L_{\textbf{m}} \rightarrow J_{0}(\n)_{\mathbf{m}} \rightarrow J_{0}(\n),
\end{equation}
where $L_{\mathbf{m}}$ is defined by the quotient
\[
 0 \rightarrow \mathbb{G}_{m} \rightarrow \prod_{i=0}^{r} \text{Res}_{K(P_i)/K} \mathbb{G}_{m}  \rightarrow L_{\mathbf{m}} \rightarrow 0.
\]
After applying the exact functor $H^{1}(-, K)$, and the left-exact functor $\text{Hom}(\Q/\Z, -)$ to the sequence in Equation~\eqref{eq:generalisedJacobian}, we obtain 
the exact sequence
\begin{equation} \label{eq:torsiongeneralisedJacobian} 0 \rightarrow \bigoplus_{i=0}^{r-1} K(P_{i})^{\times}_{\text{tors}} \rightarrow J_{0}(\n)_{\mathbf{m}}(K)_{\text{tors}} \rightarrow J_{0}(\n)(K)_{\text{tors}}   \xrightarrow{\delta}     \bigoplus_{i=0}^{r-1} K(P_{i})_{\text{tors}} \otimes \Q/\Z.\end{equation}
By Theorem \ref{thm:field} and  \cite[Proposition 5.2]{Hayes}, we have 
$K(P_{i})^{\times}_{\text{tors}} \simeq \F_{q}^{\times} \simeq \Z/(q-1)\Z$. Hence, to understand $J_{0}(\n)_{\mathbf{m}}(K)_{\text{tors}}$, we need to understand the kernel of the map $\delta$ appearing in Equation~\eqref{eq:torsiongeneralisedJacobian}. 
The key lemma to do so is Lemma 2.3.1 in \cite{yamazaki2016rational}. We present their result here in the form that will be useful to us. 

\begin{lemma}[{cf. \cite[Lemma~2.3.1]{yamazaki2016rational}}] \label{lem:delta} Fix $\n=\p^{r}$ with $\p \subseteq \mathcal{A}$ monic irreducible and $r$ a positive integer. Let $D = \sum_{i=0}^{r} a_{i}\cdot P_{i} \in \textup{Div}^{0}(X_{0}(\n))$ be a degree-zero divisor supported on $\textbf{m}$ such that $[D] \in J_{0}(\n)(K)_{\tor}$. Let $n \in \Z_{>0}$ so that $n\cdot[D]= 0$, i.e., so that there exists $f \in K(X_{0}(\n))$ such that $\Div(f)=n \cdot D$. We have
\begin{equation*}
     \delta([D])= \left( \left(\frac{f}{t_{P_{r}}^{na_{r}}}\right)(P_{r})\left(\frac{t_{P_{i}}^{na_{i}}}{f}\right)(P_{i}) \right)_{i=0}^{r-1} \otimes \frac{1}{n} \in \bigoplus_{i=0}^{r-1} K(P_{i})^{\times} \otimes \Q/\Z, 
\end{equation*}
 where $t_{P_{i}}$ is a uniformiser at the closed point $P_{i} \in X_{0}(\n)$. Notice that this description does not depend on the choice of $f$ nor on that of $t_{P_{i}}$. 
\end{lemma}
Since the map $\delta$ is a group homomorphism, we have the decomposition \begin{equation} \label{eq:kerhomomorphism} 
\ker(\delta|_{J_{0}(\n)(\Q)_{\tor}}) = \bigoplus_{\ell \text{ prime}} \ker(\delta|_{J_{0}(\n)(\Q)_{\tor}})[\ell^{\infty}] = \bigoplus_{\ell \text{ prime}} \ker(\delta|_{J_{0}(\n)(\Q)_{\tor}[\ell^{\infty}]}). \end{equation} 
By Theorem~\ref{thm:kevintorsion}, for all $\ell \nmid q(q-1)$ we have 
\begin{equation} \label{eq:kerhomomorphism} \ker(\delta|_{J_{0}(\n)(\Q)_{\tor}[\ell^{\infty}]}) =  \ker(\delta|_{C(\n)[\ell^{\infty}]})=\ker(\delta|_{C(\n)})[\ell^{\infty}].\end{equation}
In other words, to compute $\ker(\delta|_{J_{0}(N)(\Q)_{\tor}[\ell^{\infty}]})$, it is enough to understand the map $\delta$ restricted to the subgroup $C(\n)$. The next proposition takes a closer look at $\delta$ when applied to elements of $C(\n)$ when $\n=\p^r$ is a prime power. 
\begin{proposition} \label{prop:coef} Fix $\n=\p^{r}$ with $\p \subseteq \mathcal{A}$ monic irreducible and $r$ a positive integer. Let $D \in \mathcal{C}(\n)$ be a cuspidal divisor such that for some integer $n(D) \in \Z_{>0}$ we have $g(n(D) \cdot D)= \prod_{\m \mid \n} \Delta_{\m}(z)^{k_{\m}} \in K(X_{0}(\n))$ for some $k_{\m} \in \Z$. We have 
$$ \delta(D) = \left( \prod_{\m \mid \n}\left(\frac{\gcd(\p^{k}, \m)}{\m}\right)^{-(q^{2}-1)k_{\m}}\prod_{\m \mid \n }e_{L}\left(\frac{\eta_{\m} \gcd(\p^{k}, \m)}{\m}\right)^{-k_{\m}} \right)_{k=0}^{r-1} \otimes \frac{1}{n(D)},$$
where $\eta_{\m}\in \mathcal{A}$ and  $e_{L}$ is the Carlitz exponential associated to a lattice $L$ as defined in ~\cite{GekelerProduct}. 
\end{proposition}

\begin{proof}
From the explicit expression in Lemma~\ref{lem:delta} we see that to prove the result we want to find a formula for $\left(\frac{g(n(D)\cdot D)(z)}{t_{P_{i}}^{n(D)a_{i}}}\right)(P_{i})$ for all closed points $P_{i}$. Since $P_{i}$ is a Galois orbit, it is enough to consider one representative $\omega\in\cusp(X_{0}(\n)(\C_{\infty}))$ for each $i$.  Take a cusp representative $\omega =\left[\begin{smallmatrix}
    \fa' \\ \m'
\end{smallmatrix}\right] \in \cusp(X_{0}(\n)(\C_{\infty}))$ and a matrix $A \in \SL_{2}(\mathcal{A})$ such that $A^{-1} \cdot \omega_{\infty} = \left[\begin{smallmatrix}
    \fa' \\ \m'
\end{smallmatrix}\right]$. The Fourier expansion of $\Delta_{\m}(z)$ at $\omega$ is the expansion of $\Delta_{\m}(A^{-1}z)$ at $\omega_{\infty}$.  Take $t_{\omega_{\infty}} \in K (X_{0}(\n))^{\times}$ to be a
uniformiser at $\omega_{\infty}$ which has 1 as the leading term of its Fourier expansion.
Let $h(z) := \prod_{\m \mid \n} \Delta_{\m}(z)^{k_{\m}} \in \mathcal{E}_{\n} \otimes_{\Z}\Q$ be any $\Delta$-quotient. Then we have
$$ \left(\frac{h(z)}{t_{\omega_{\infty}}^{\text{ord}_{\infty}(g)}}\right)(\omega_{\infty})= 1.$$
Now take the uniformiser at $[\omega]$ given by
$$ t_{\omega} := t_{\omega_{\infty}} \circ A \in K(X_{0}(\n)).$$
The identity in  Lemma~\ref{lem:identityDelta} yields
\begin{align*}
&\left(\frac{h(z)}{t_{\omega}^{\text{ord}_{\omega}(g)}}\right)(\omega)=\left(\frac{h(z)\circ A^{-1}}{t_{\omega_{\infty}}^{\text{ord}_{\omega}(g)}}\right)(\omega_{\infty})= \\ &\prod_{\m \mid \n}\left(\frac{\gcd(\p^{k}, \m)}{\m}\right)^{(q^{2}-1)k_{\m}} \frac{\prod_{\m \mid \n }\Delta\left(\frac{\gcd(\p^{k}, \m)^{2}}{\m}z+\frac{\eta_{(\m, \p^{k})} \gcd(\p^{k}, \m)}{\m}\right)^{k_{\m}}}{t_{\omega_{\infty}}^{\text{ord}_{\omega}(g)}}(\omega_{\infty});\end{align*}
where $\eta_{(\m, \p^{k})}$ is some element in $ \mathcal{A}$.
The product formula for the function $\Delta$ found by Gekeler in \cite[Theorem pg.137]{GekelerProduct} finishes the proof. 
\end{proof}

\begin{remark} \label{rmk:subgroupimage}
Take the notation as in the statement of Proposition~\ref{prop:coef}. Note that the Fourier expansion of $\Delta$ at a cusp of height $\m=p^{i}$ has $K(P_{i})$-rational coefficients, see \cite{GekelerProduct}. For a fixed $\m \mid \n$, recall that $\fd=\fd(\m)=\gcd(\m, \n/\m)$. Let $\mathcal{O}_{\fd}\coloneqq \langle \p \rangle \subseteq (K(\mathfrak{d})^{+})^{\times}$ be the subgroup of $(K(\fd)^{+})^{\times}$ generated by $\p$.
The previous proposition shows that $\im(\delta|_{\mathcal{C}(\n)}) \subseteq \bigoplus_{i=0}^{r-1} \left(\mathcal{O}_{\fd}\cdot \mu(K(P_{i}))\right) \otimes \Q/\Z $. 
\end{remark}

For the rest of the paper we fix $\n=\p^{r}$ with $\p \subseteq \mathcal{A}$ monic irreducible and $r$ a positive integer.

\begin{definition} Consider the subgroup $\Lambda_{\mathfrak{d}}$ given by the image of $\mathcal{O}_{\mathfrak{d}}$ in the quotient $K(P_{i})^{\times} / \mu(K(P_{i}))$. We denote $\Lambda := \bigoplus_{\m \mid \n} \Lambda_{\fd(\m)}$. 
\end{definition}

By Remark~\ref{rmk:subgroupimage} and since all torsion units $x \in \mu(K(P_{i}))$ with $n \mid \ord(x)$ satisfy $x \otimes \frac{a}{b} = 1 \otimes \frac{a}{nb}$ and the latter is in turn the identity element in $K(P_{i}) \otimes \Q/\Z$, the map $\delta$ factorises as
\begin{equation} \label{fig:fact}
\begin{tikzcd}  
C(\n)  {} \arrow[rrr, "\overline{\delta}"] \arrow[rrrdd, "\delta" '] &  &  & {}   \Lambda\otimes \Q/\Z     \arrow[dd, "\phi"] \\
      &                                &  &                      \\
      &                                &  & {} \bigoplus_{i=0}^{r-1} K(P_{i})^{\times} \otimes \Q/\Z         
\end{tikzcd}
\end{equation}

\begin{lemma} \label{lem:torsunits}
    The map $\phi$ in the diagram of ~\eqref{fig:fact} is an injection. 
\end{lemma}

\begin{proof}
For each $\m=p^{i}$ with $\m \mid \n$, let $\fd= \fd(\m) = \gcd(\m, \n/\m)$ as before. Consider the following tensoring map on $\Z$-modules:
\[ \phi_{\fd(\m)} : \Q/\Z \longrightarrow K(P_{i})^{\times} \otimes \Q/\Z, \quad x \mapsto \p \otimes x.\]
If we prove that for each $\fd$ the map $\phi_{\fd}$ is a split injection, then the injectivity of $\phi$ follows since $\Lambda \otimes \Q/\Z \simeq \bigoplus_{\fd} (\Lambda_{\fd} \otimes \Q/\Z)$ and $\left(\bigoplus_{i} K(P_{i}^{\times}) \right)\otimes \Q/\Z \simeq \bigoplus_{i} (K(P_{i})^{\times} \otimes \Q/\Z)$. 
Since $\Q/\Z$ is an injective module, splitting of $\phi_{\fd}$ follows automatically if we prove its injectivity.
Now the prime $\p \in K(P_i)$
is not torsion, thus its image $\mathcal{O}_{\fd}$ is a free $\Z$-module (note that $\F_{q}$ is the only torsion part in $K(P_i)^{\times}$). Hence, the $\phi_{\fd}$ is injective for all $\m$. 
\end{proof}

\begin{remark}
    The factorisation given in Equation~\eqref{fig:fact} together with Lemma~\ref{lem:torsunits} show that $\ker(\delta|_{C(\n)}) = \ker(\overline{\delta}|_{C(\n)})$. Since the map $\overline{\delta}$ no longer takes the torsion units into account, it is relatively easier to compute. In the next section we will explicitly compute the images of  the generators in Theorem~\ref{thm:cuspgen} under $\overline{\delta}$.
\end{remark}

As a result of a technical study of $\overline{\delta}$ we obtain the main result which we recall below. 
\begin{theorem} Fix $\n=\p^{r}$ with $\p \subseteq \mathcal{A}$ monic irreducible and $r$ a positive integer. The map $\overline{\delta}|_{\CC(\n)}$ is an injection. In particular, for all odd primes $\ell$ with $\ell \nmid q(q^{2}-1)$ the subgroup $J_{0}(\p^{r})_{\mathbf{m}}(\Q)[\ell^{\infty}]$ is trivial. 
\end{theorem}

\begin{proof}
    This is a consequence of Equations \eqref{eq:kerhomomorphism}, Lemma~\ref{lem:torsunits} and Theorem~\ref{prop:maximalrank} proven in Section~\ref{sec:kerneldelta}, which states that $\overline{\delta}|_{\CC(\n)}$ is injective.
\end{proof}

The goal of the coming two sections is to prove Theorem~\ref{prop:maximalrank}. In Section~\ref{sec:imagedelta} we compute the image of the cuspidal divisors appearing in Theorem~\ref{thm:cuspgen}. This will allow us to find a matrix representation of the map $\overline{\delta}|_{\CC(\p^r)}$. In Section~\ref{sec:kerneldelta} we will use this representation to determine the kernel of $\overline{\delta}|_{\CC(\n)}$. 

\section{The image of $\overline{\delta}$.}
\label{sec:imagedelta}
Now we proceed with the computations of the map $\overline{\delta}$. The main goal of this section is to compute $\im(\overline{\delta}|_{\CC(\p^r)})$ in Proposition~\ref{prop:imageCpr}. 
\begin{remark}\label{rmk:torsionunits}
Recall from Proposition~\ref{prop:coef} that the product $\prod_{\m \mid \n}\left(\frac{\gcd(\p^{k}, \m)}{\m}\right)^{(q^{2}-1)k_{\m}} \prod_{\m \mid \n }e_{L}\left(\frac{\eta_{\m} \gcd(\p^{k}, \m)}{\m}\right)^{k_{\m}}$ appears in the expression of $\delta(\overline{D})$
 for $D \in \CC(\n)$. From \cite[Theorem 8.4]{conradcarlitz} we have that the factor of the product $\prod_{\m \mid \n }e_{L}\left(\frac{\eta_{\m} \gcd(\p^{k}, \m)}{\m}\right)^{k_{\m}} \in K(P_{i})$ is a torsion unit. Hence, by Lemma~\ref{lem:torsunits}, in order to show that the map $\overline{\delta}$ in Equation~\eqref{fig:fact} is injective we will only need to control the product $ \prod_{\m \mid \n}\left(\frac{\gcd(\p^{k}, \m)}{\m}\right)^{(q^{2}-1)k_{\m}}$ in the product expression. 
\end{remark}

The product expression appearing Proposition \ref{prop:coef} together with Remark~\ref{rmk:torsionunits} is used in the following for the computation of the images of the divisors $\{C_{i}\}_{i=1}^{\lfloor\frac{r-1}{r}\rfloor}$,  $\{C_{i}-C_{i-1}\}_{i=\lfloor\frac{r+1}{2}\rfloor}^{r-2}$ and $\{D_{i}\}_{i=0, r-1}$ under the map $\overline{\delta}$. 
Letting $r \geq 2$, by Theorem~\ref{thm:cuspgen} we have 
$$\CC(\p^r) = \left(\bigoplus_{i=1}^{m}\langle \overline{C_i}\rangle\right)\oplus\left(\bigoplus_{i=m+1}^{r-2}\langle \overline{C_i-\pp C_{i+1}}\rangle\right)\oplus \langle \overline{D_{r-1}}\rangle\oplus \langle  \overline{D_0}\rangle,$$ where $m:=\lfloor \frac{r-1}{2}\rfloor$.

\begin{notation}
Let $\sigma \in \Lambda$ a tuple of the form $$ \sigma = (\p^{\sigma_{k}})_{k=0}^{r-1};$$
where each for each $k$, the element $\p^{\sigma_{k}}$ lies in $\Lambda_{\fd}$ with $\fd=\fd(\p^{k})$.
\end{notation}

\begin{proposition} \label{prop:imageCpr} The images of $\overline{C_{i}}$ for $i \in \{ 1, \ldots, \lfloor\frac{r-1}{2}\rfloor\}$ under the map $\overline{\delta}$ are given by
\begin{equation}\label{eq:deltasmall}
    \overline{\delta}(\overline{C_{i}}) = (\p^{-\sigma(i)_{k}} )_{k=0}^{r-1} \otimes \frac{q+1}{\pp^{r-i}(\pp^{2}-1)},
\end{equation}
where 
\begin{equation}\label{eq:deltasmall}
    \sigma(i)_{k} = \begin{cases}
    (r-i)(\pp^2-2\pp+1)+\pp-1 & \text{if } 0 \leq k \leq i-1; \\
    (r-i)(\pp^2-2\pp+1)+2\pp-1 & \text{if } k=i;\\
    (r-k)(\pp^2-2\pp+1)+\pp-1 & \text{if } i+1 \leq k \leq r-1. \\
    \end{cases}
\end{equation}
Similarly, for $i \in \{\lfloor\frac{r+1}{2}\rfloor, \ldots, r-2\}$ we have
\begin{equation} \label{eq:deltabig}
\overline{\delta}(\overline{C_i-\pp C_{i+1}}) = (\p^{-\sigma(i)_{k}})_{k=0}^{r-1} \otimes \frac{q+1}{\pp^{i}(\pp^{2}-1)}, \end{equation}
where
\begin{equation} \label{eq:deltabig}
\sigma(i)_{k}= \begin{cases}
    (\pp-1)^{2} & \text{if } 0 \leq k\leq i-1; \\
    \pp^2-\pp+1 & \text{if }k=i; \\
    -\pp & \text{if } k=i+1; \\
    0 & \text{if } i+2 \leq k\leq r-1. \\
\end{cases}\end{equation}
Finally, we have
\begin{equation} \label{eq:delta 0}
\overline{\delta}(\overline{D_{0}})=  (\p^{-(q-1)\sigma(0)_{k}})_{k=0}^{r-1} \otimes \frac{q+1}{\pp-1}, \quad \text{ where } \quad \sigma(0)_{k} = 1  
\end{equation}
for all $0 \leq k \leq r-1$,
and
\begin{equation} \label{eq:delta r-1} \overline{\delta}(\overline{D_{r-1}}) = (\p^{-\sigma(r-1)_{k}})_{k=1}^{r-1} \otimes \frac{q+1}{\pp^{2}-1},\end{equation}
where the exponents $\sigma(r-1)_{k}$ can be found in Appendix
\ref{sec:appendix} -- see Remark~\ref{rmk:sigma(r-1)}. 
\end{proposition}

\begin{proof}
We first determine $\overline{\delta}(\overline{C_{i}})$ for $i \in \{1, \ldots, \lfloor\frac{r-1}{2}\rfloor\}$. Let $m=\min\{i, r-i\}$. Then \cite[Theorem 3.6]{kevincusp} states that
\[g(C_i) = \left(\frac{\Delta_{\p^{i}}^{\pp^{2}+1}\Delta_{\p^{r-1}}^{\pp-1}\Delta_{\p^{r}}^{\pp}}{\Delta_{\p^{i-1}}^{\pp}\Delta_{\p^{i+1}}^{\pp}\Delta_{\p^{r}}^{\pp^{2}}}\right)^{\frac{1}{\pp^{r-m}(\pp^{2}-1)(q-1)}}.\]
Since by Theorem~\ref{thm:cuspgen} we have 
$\ord(\overline{C_i})=\pp^{r-i}M(\p)$, using Proposition~\ref{prop:coef} 
we have 
\begin{equation*}
    \overline{\delta}(\overline{C_{i}}) = (\p^{-\frac{\sigma(i)_{k}}{q-1}} )_{k=0}^{r} \otimes \frac{q^2-1}{\pp^{r-i}(\pp^{2}-1)}=(\p^{-\sigma(i)_{k}} )_{k=0}^{r} \otimes \frac{q+1}{\pp^{r-i}(\pp^{2}-1)}
\end{equation*}
with $\sigma(i)_{k}$ as in  Equation~\eqref{eq:deltasmall}.

On the other hand, for $i=\lfloor\frac{r+1}{2}\rfloor, \ldots, r-2$ and $m=\min\{i, r-i\}$, \cite[Lemma 3.8]{kevincusp} computes
\[g(C_i-\pp C_{i+1}) = \left(\frac{\Delta_{\p^{i}}^{\pp^{2}+\pp+1}\Delta_{\p^{i+2}}^{\pp}}{\Delta_{\p^{i-1}}^{\pp}\Delta_{\p^{i+1}}^{\pp^{2}+\pp+1}}\right)^{\frac{1}{\pp^{i}(\pp^{2}-1)(q-1)}}.\]
Since by Theorem~\ref{thm:cuspgen} we have
$\ord(\overline{C_i-\pp C_{i+1}})=\pp^{i}M(\p)$,
using Proposition~\ref{prop:coef} 
we obtain
\begin{equation*}
\overline{\delta}(\overline{C_i-\pp C_{i+1}}) = (\p^{-\frac{\sigma(i)_{k}}{q-1}})_{k=0}^{r} \otimes \frac{q^2-1}{\pp^{i}(\pp^{2}-1)}=(\p^{-\sigma(i)_{k}})_{k=0}^{r} \otimes \frac{q+1}{\pp^{i}(\pp^{2}-1)}
\end{equation*}
with $\sigma(i)_{k}$ as in Equation~\eqref{eq:deltabig}.

We are left with computing $\overline{\delta}(\overline{D_{i}})$ for $i=0, r-1$. We start with $i=0$. By 
\cite[Lemma 3.9]{kevincusp}, we have 
$$ g(D_{0}) = \left(\frac{\Delta_{\p^{r-1}}}{\Delta_{\p^r}}\right)^{\frac{1}{\pp-1}},$$
and Theorem~\ref{thm:cuspgen} gives $\ord(\overline{D_{0}}) = N(\p)$. Hence, Proposition~\ref{prop:coef} yields Equation~\eqref{eq:delta 0}.
Finally for $i = r-1$ we have a case distinction for the definition of $D_{r-1}$ -- see \cite[Section 3.2]{kevincusp}. The modular units $g(D_{r-1})$ are given in \cite[Lemma 3.8]{kevincusp} for all different cases. 
Using this lemma, if $\sigma(r-1)_{k}$ denotes the exponent at $\p$ of the leading coefficient of the $q$-expansion of the modular unit $g(M(\p)\cdot D_{r-1})$ at the cusp $\frac{1}{\p^{k}}$ for $0 \leq k \leq r-1$, then we have
$$ \overline{\delta}(\overline{D_{r-1}}) = (\p^{
-\sigma(r-1)_{k}})_{k=0}^{r-1} \otimes \frac{q+1}{\pp^{2}-1}.$$
For convenience of the reader, the computation of the entries $\sigma(r-1)_{k}$ for $0 \leq k \leq r-1$ can be found in the Appendix \ref{sec:appendix}. 
\end{proof}


\section{Injectivity of the map $\delta$.}
\label{sec:kerneldelta}

In this section, we show that the map $\overline{\delta}$ is injective to complete the proof of Theorem~\ref{thm:mainthm}. Using Proposition~\ref{prop:imageCpr}, we will first find a matrix representing the map $\overline{\delta}$. To do so, we introduce the following notation. 

\begin{notation} 
Recall that for a divisor $D \in \{C_{i}: 1 \leq i \leq \lfloor{\frac{r-1}{2}\rfloor}\}\cup\{C_{i}-\pp C_{i+1}:\lfloor{\frac{r+1}{2}\rfloor} \leq i \leq r-2\}\cup\{D_{r-1}\}$, we denote by $\sigma(i)_{k}$ the exponent of $\p$ at the $k$-th coordinate of the tuple of element in the tensor product $\oplus_{k=0}^{r-1}K(P_{k})^{\times} \otimes \Q/\Z$ of the representatives $\overline{\delta}(\overline{D})$ given in Equations~\eqref{eq:deltasmall}, \eqref{eq:deltabig} and \eqref{eq:delta r-1}.
\end{notation}

Notice that, since we can take $\frac{q+1}{\pp^{r}(\pp^2-1)} \in \Q/\Z$ as common factor of the images of the divisors $\{C_{i}\}_{i=1}^{\lfloor\frac{r-1}{2}\rfloor}$, $\{C_{i}-\pp C_{i+1}\}_{i=\lfloor\frac{r+1}{2}\rfloor}^{r-2}$ and $D_{r-1}$ computed in Proposition~\ref{prop:imageCpr}, the $r \times r$ matrix $(M_{\delta}(i, j))_{i, j}$ given by \[
\mathbf{M}_{\delta}(i, j) := -1\cdot\begin{cases}
    \pp^{r}(q-1)\sigma(0)_{j} & \text{for } i=0,\ j=0, \ldots, r-1; \\
    \pp^{i}\sigma(i)_{j} & \text{for } i=1, \ldots, \lfloor\frac{r-1}{2}\rfloor,\ j=0, \ldots, r-1 ;\\
    \sigma(i)_{j} & \text{for } i=\lfloor\frac{r+1}{2}\rfloor, \ldots, r-1,\  j=0, \ldots, r-1;\\
    \pp^{r} \sigma(i)_{j} & \text{for } i=r-1,\ j=0, \ldots, r-1;\\
\end{cases}\]
represents the map $\overline{\delta}$ defined in Equation~\eqref{fig:fact}. Notice that the entries of $\mathbf{M}_{\delta}$ are labelled from $0$ to $r-1$ instead of the ususal $1$ to $r$. To finish the proof of Theorem~\ref{thm:mainthm}, we need to check the injectivity of $\overline{\delta}$. That is, we want to show that $\text{det}(\mathbf{M}_{\delta}) \neq 0$.  

\begin{notation}
    For a matrix $M=(m_{i, j})_{0 \leq i \leq n-1, 0 \leq j \leq m-1}$ of dimension $n\times m$, we use the following notation throughout the section. 
    \begin{itemize}
        \item For $0 \leq i \leq n-1$, we denote by $M[i]$ the $i$-th row of $M$.
         \item For $0 \leq j \leq m-1$, we denote by $M[[j]]$ the $j$-th column of $M$.
          \item For $0 \leq i \leq n-1$ and $0 \leq j \leq m-1$, we denote by $M(i,j)$ the $(i, j)$-th entry of $M$.
          \item For a given $i$, we represent the row $M[i]$ of $M$ as $M[i]=(\ast_{0}, \dots, \ast_{m-1})$, where $\ast_{j}$ is the entry $M(i, j)$. We will write the subindex $j$ only when necessary to clarify the position of the entries of a given tuple. 
    \end{itemize}
\end{notation}

\begin{theorem} \label{prop:maximalrank}
    The matrix $\mathbf{M}_{\delta}$ representing the map $\delta$ has maximal rank. 
\end{theorem}


We will take the rest of the section to prove Theorem~\ref{prop:maximalrank}. To do so, ill use the following auxiliary matrix. Let 
$$M_{\delta}(i, j) := \begin{cases}
    \sigma(0)_{j} & \text{for } i=0, \ j=0, \ldots, r-1; \\
    \sigma(i)_{j} & \text{for } i=1, \ldots, \lfloor\frac{r-1}{2}\rfloor,\  j=0, \ldots, r-1 ;\\
    \sigma(i)_{j} & \text{for } i=\lfloor\frac{r+1}{2}\rfloor, \ldots, r-1, \  j=0, \ldots, r-1;\\
    \sigma(r-1)_{j} & \text{for } i=r-1, \ j=0, \ldots, r-1.\\
\end{cases}$$
The matrices $\mathbf{M}_{\delta}$ and $M_{\delta}$ differ only from each other by rescaling the rows. Hence, to check that $\mathbf{M}_{\delta}$ has non-zero determinant, it is enough to check that $\text{det}(M_{\delta}) \neq 0$.
If $r \leq 6$, showing $\text{det}(M_{\delta}) \neq 0$ can be done by a straightforward computation. However, to prove the result in more generality for $r \geq 7$, we will first transform the matrix $M_{\delta}$ into a matrix $M_{\delta}^{h}$ with the same determinant that is easier to manipulate. We will do that in two steps: we will first do a row reduction to produce many zeroes above the main diagonal of the matrix, and then we will do a column reduction to create many zeroes in the first column. The combination of these steps will be enough to allow us to prove that $\det(M_{\delta}) \neq 0$. 

\textbf{Step 1:} Let $r \geq 7$. Let $M^{H}_{\delta}$ be the matrix obtained from $M_{\delta}$ after applying the following \underline{row} transformations:
\begin{itemize}
    \item $M^{H}_{\delta}[1] = M_{\delta}[1]- M_{\delta}[2]$, for $i=1$;
    \item $M^{H}_{\delta}[i] = 2M_{\delta}[i]-M_{\delta}[i-1] - M_{\delta}[i+1]$, for $1 < i \leq \lfloor\frac{r-1}{2}\rfloor-1$;
    \item $M^{H}_{\delta}[\lfloor\frac{r-1}{2}\rfloor] = 2M_{\delta}[\lfloor\frac{r-1}{2}\rfloor]+\sum_{k=1}^{\lfloor r-3/2\rfloor} (p^k-1)M_{\delta}[r-1-k] - M_{\delta}[(r-3)/2]$;
    \item $M^{H}_{\delta}[\lfloor\frac{r+1}{2}\rfloor] = M_{\delta}[\lfloor\frac{r+1}{2}\rfloor]-\left(M_{\delta}[\lfloor\frac{r-1}{2}\rfloor]+\sum_{k=1}^{\lfloor r-3/2\rfloor} (p^k-1)M_{\delta}[r-1-k]\right)$;
    \item $M^{H}_{\delta}[i] = M_{\delta}[i]-M_{\delta}[i-1]$, for $ \lfloor\frac{r+1}{2}\rfloor < i \leq r-1$.
\end{itemize}

\begin{claim} \label{claim:matrix1} The matrix $M^{H}_{\delta}$ is a matrix 
of the form 

$$\begin{psmallmatrix}
1 & 1 & 1 & 1 & 1 & \dots & \dots & \dots & \dots& \dots & \dots& 1 \\
 (\pp-1)^2 & \pp^2-\pp+1 & -\pp & 0 & 0 & \dots & \dots & \dots & \dots& \dots & \dots& 0 \\
 0& -\pp & \pp^2+1 & -\pp & 0 & \dots & \dots & \dots &\dots & \dots& \dots& 0 \\
\dots & \dots & \dots & \dots &\dots & 0 & \dots  & \dots & \dots& \dots&\dots& 0 \\
 0 & 0 & \dots & -\pp & \pp^2+1 & -\pp & 0 & \dots & \dots& \dots&\dots& 0 \\
\alpha & \dots & \dots & \alpha & \alpha+\pp & \beta+2\pp & \pp \left(\pp^{\lfloor \frac{r-1}{2}\rfloor}-1\right) & 0 & \dots& \dots&\dots& 0 \\
 -\alpha & -\alpha & \dots & \dots& -\alpha & -\alpha-\pp & -\pp^{\lfloor\frac{r+1}{2}\rfloor}+\pp^2+1 & -\pp & 0 & \dots &\dots &0\\
  \dots & \dots & \dots & \dots &\dots & \dots & \dots  & \dots & \dots& \dots& \dots &  0 \\
   0 & \dots  & \dots &\dots &\dots &\dots &\dots &\dots & -\pp & \pp^2+1 & -\pp &0  \\
 0 & \dots & \dots & \dots &\dots &\dots &\dots &\dots &\dots & -\pp & \pp^2+1 & -\pp \\
 \sigma(r-1)_{0} & \sigma(r-1)_{1} & \dots & \dots & \dots & \dots & \dots & \dots & \dots & \sigma(r-1)_{r-3}& \sigma(r-1)_{r-2} &\sigma(r-1)_{r-1}\\
\end{psmallmatrix}$$ 
where $\alpha= \pp^{\lfloor \frac{r+1}{2}\rfloor}-\pp^{\lfloor \frac{r-1}{2}\rfloor}$ and $\beta= \pp^{\lfloor \frac{r+1}{2}\rfloor}-\pp^{\lfloor \frac{r-1}{2}\rfloor}+\pp^2-2 \pp+1.$
More formally, we have 
\begin{center}
    $M_{\delta}^{H}[i] = \begin{cases}
    (1, \ldots, 1) & \text{if } i=0; \\
    ((\pp-1)^{2}, \pp^{2}-\pp+1, -\pp, 0, \ldots, 0) & \text{if } i=1; \\
    (0, \ldots, 0, -\pp_{i-1}, (\pp^{2}+1)_{i}, -\pp_{i+1}, 0, \ldots, 0) & \text{if } 1 < i < \lfloor\frac{r-1}{2}\rfloor; \\
    (\alpha, \ldots, \alpha, (\alpha+\pp)_{\lfloor\frac{r-3}{2}\rfloor}, \beta_{\lfloor\frac{r-1}{2}\rfloor}, \pp \left(\pp^{\lfloor \frac{r-1}{2}\rfloor}-1\right)_{\lfloor\frac{r+1}{2}\rfloor}, 0, \ldots, 0) & \text{if } i=\lfloor\frac{r-1}{2}\rfloor; \\
    (-\alpha, \ldots ,-\alpha, (-\alpha-\pp)_{\lfloor\frac{r-1}{2}\rfloor}, (-\pp^{\lfloor\frac{r+1}{2}\rfloor}+\pp^{2}+1)_{\lfloor\frac{r+1}{2}\rfloor}, -\pp_{\lfloor\frac{r+3}{2}\rfloor}, 0, \ldots, 0) & \text{if } i=\lfloor\frac{r+1}{2}\rfloor; \\
    (0, \ldots, 0, -\pp_{i-1}, (\pp^{2}+1)_{i}, -\pp_{i+1}, 0, \ldots, 0) & \text{if } \lfloor\frac{r+1}{2}\rfloor < i < r-1; \\
    (\sigma(r-1)_{j})_{j=0}^{r-1} & \text{if } i=r-1.
    \end{cases}$
\end{center}
\end{claim}

\begin{proof} Notice that, for $1 < i \leq \lfloor\frac{r-1}{2}\rfloor-1$, we can rewrite the transformations as follows: 
$$ M^{H}_{\delta}[i] = (M_{\delta}[i]-M_{\delta}[i+1]) - (M_{\delta}[i-1] - M[i]).$$
Hence, it is enough to look at the differences $M_{\delta}[i] - M_{\delta}[i+1]$ for $1 \leq i \leq \lfloor\frac{r-1}{2}-1\rfloor$. 
From Equation~\eqref{eq:deltasmall} it follows that
$$M_{\delta}[i] - M_{\delta}[i+1] =
((\pp^2-2\pp+1), \ldots, (\pp^2-2\pp+1)_{i-1}, (\pp^2-\pp+1)_{i}, -\pp_{i+1}, 0, \ldots, 0),$$
for $1 \leq i \leq \lfloor\frac{r-3}{2}\rfloor$. Hence, we have 
$$M_{\delta}^{H}[1]=((\pp-1)^{2}, \pp^{2}-\pp+1, -\pp, 0, \ldots, 0)$$
for $i=1$ and 
$$(M_{\delta}[i]-M_{\delta}[i+1]) - (M_{\delta}[i-1] - M[i])=(0, \ldots, 0, -\pp_{i-1}, (\pp^{2}+1)_{i}, -\pp_{i+1}, 0, \ldots, 0)$$
for $1<i< \lfloor\frac{r-3}{2}\rfloor$.
Similarly, from Equation~\eqref{eq:deltabig} it follows that
$$ M^{H}_{\delta}[i] = M_{\delta}[i]-M_{\delta}[i-1] = (0, \ldots, 0, -\pp_{i-1}, (\pp^{2}+1)_{i}, -\pp_{i+1}, 0, \ldots, 0),$$ 
for $ \lfloor\frac{r+1}{2}\rfloor < i \leq r-1$. Hence, we are left with computing $M_{\delta}^{H}[\lfloor\frac{r-1}{2}\rfloor]$ and $M_{\delta}^{H}[\lfloor\frac{r+1}{2}\rfloor]$. Recall from \textbf{Step 1} that 
    $$M^{H}_{\delta}[\lfloor\frac{r-1}{2}\rfloor] = 2M_{\delta}[\lfloor\frac{r-1}{2}\rfloor]+\sum_{k=1}^{\lfloor r-3/2\rfloor} (p^k-1)M_{\delta}[r-1-k] - M_{\delta}[(r-3)/2]=$$
    $$\left(M_{\delta}[\lfloor\frac{r-1}{2}\rfloor]+\sum_{k=1}^{\lfloor r-3/2\rfloor} (p^k-1)M_{\delta}[r-1-k]\right) - \left(M_{\delta}[(r-3)/2]-M_{\delta}[\lfloor\frac{r-1}{2}\rfloor]\right) $$
    and
    $$M^{H}_{\delta}[\lfloor\frac{r+1}{2}\rfloor] = M_{\delta}[\lfloor\frac{r+1}{2}\rfloor]-\left(M_{\delta}[\lfloor\frac{r-1}{2}\rfloor]+\sum_{k=1}^{\lfloor r-3/2\rfloor} (p^k-1)M_{\delta}[r-1-k]\right).$$
Hence if we show that
\begin{equation} \label{eq:equationclaim1} M_{\delta}[\lfloor\frac{r-1}{2}\rfloor]+\sum_{k=1}^{\lfloor r-3/2\rfloor} (p^k-1)M_{\delta}[r-1-k] =  ( \beta, \dots, \beta, (\beta+\pp)_{\lfloor\frac{r-1}{2}\rfloor},  \pp \left(\pp^{\lfloor\frac{r-1}{2}\rfloor}-1\right)_{\lfloor\frac{r+1}{2}\rfloor}, 0 , \dots, 0 )\end{equation}
then the equalities 
$$ M^{H}_{\delta}[\lfloor\frac{r-1}{2}\rfloor] = (\alpha, \ldots, \alpha, (\alpha+\pp)_{\lfloor\frac{r-3}{2}\rfloor}, (\beta+2\pp)_{\lfloor\frac{r-1}{2}\rfloor}, \pp \left(\pp^{\lfloor \frac{r-1}{2}\rfloor}-1\right)_{\lfloor\frac{r+1}{2}\rfloor}, 0, \ldots, 0)$$ 
and 
$$M^{H}_{\delta}[\lfloor\frac{r+1}{2}\rfloor] =(-\alpha, \ldots ,-\alpha, (-\alpha-\pp)_{\lfloor\frac{r-1}{2}\rfloor}, (-\pp^{\lfloor \frac{r+1}{2}\rfloor}+\pp^{2}+1)_{\lfloor\frac{r+1}{2}\rfloor}, -\pp_{\lfloor\frac{r+3}{2}\rfloor}, 0, \ldots, 0)$$
follow and finish the proof of the claim. 
Thus, to complete the proof we only need to verify Equation ~\eqref{eq:equationclaim1}. 
Fix $j \geq \lfloor\frac{r+3}{2}\rfloor$. For all $\lfloor\frac{r+1}{2}\rfloor \leq i \leq r-2$ we have 
$M_{\delta}[i][j]=(\pp-1)^{2}$ if $i\geq j+1$, $M_{\delta}[i][j]=\pp^{2}-\pp+1$ if $i=j$, $M_{\delta}[i][j]=-\pp$ if $i=j-1$ and $M_{\delta}[i][j]=0$ if $i \leq j-2$. Hence,
$$\left(M_{\delta}[\lfloor\frac{r-1}{2}\rfloor]+\sum_{k=1}^{\lfloor r-3/2\rfloor} (p^k-1)M_{\delta}[r-1-k]\right)_{j}= $$ 
$$ (r-j)(\pp^2-2\pp+1)+\pp-1 -\pp(\pp^{r-j}-1)+(\pp^{2}-\pp+1)(\pp^{r-j-1}-1)+$$ $$(\pp-1)^{2}\sum_{k=1}^{r-2-j}(\pp^{r-1-j-k}-1)=$$
$$(r-j)(\pp^2-2\pp+1)+\pp-1 -\pp(\pp^{r-j}-1)+(\pp^{2}-\pp+1)(\pp^{r-j-1}-1)+$$ $$(\pp-1)^{2}\pp^{r-1-j}\sum_{k=1}^{r-2-j}\pp^{-k}-(\pp-1)^{2}(r-2-j)=$$
$$(r-j)(\pp^2-2\pp+1)+\pp-1 -\pp(\pp^{r-j}-1)+(\pp^{2}-\pp+1)(\pp^{r-j-1}-1)+$$ $$(\pp-1)(\pp^{r-1-j}-\pp)-(\pp-1)^{2}(r-2-j)= 0.$$
Moreover, for $j\in\{\lfloor\frac{r+1}{2}\rfloor, \lfloor\frac{r-1}{2}\rfloor\}$ we have
$$\left(M_{\delta}[\lfloor\frac{r-1}{2}\rfloor]+\sum_{k=1}^{\lfloor r-3/2\rfloor} (p^k-1)M_{\delta}[r-1-k]\right)_{\lfloor\frac{r+1}{2}\rfloor}= $$ 
$$
(r-\lfloor\frac{r+1}{2}\rfloor)(\pp^2-2\pp+1)+\pp-1+(\pp^{2}-\pp+1)(\pp^{r-\lfloor\frac{r+1}{2}\rfloor-1}-1)+$$ $$(\pp-1)^{2}\sum_{k=1}^{r-1-j}(\pp^{r-1-j-k}-1)=$$
$$ (r-\lfloor\frac{r+1}{2}\rfloor)(\pp^2-2\pp+1)+\pp-1+(\pp^{2}-\pp+1)(\pp^{r-\lfloor\frac{r+1}{2}\rfloor-1}-1)+$$ $$(\pp-1)(\pp^{r-\lfloor\frac{r+1}{2}\rfloor-1}-1)-(p-1)^2(r-1-j)=$$ $$p(p^{r-\lfloor\frac{r+1}{2}\rfloor}-1)=
p(p^{\lfloor\frac{r-1}{2}\rfloor}-1);$$
and
$$\left(M_{\delta}[\lfloor\frac{r-1}{2}\rfloor]+\sum_{k=1}^{\lfloor r-3/2\rfloor} (p^k-1)M_{\delta}[r-1-k]\right)_{\lfloor\frac{r-1}{2}\rfloor}= $$ $$
\lfloor \frac{r+1}{2} \rfloor(\pp^{2}-\pp+1)+\pp-1+(\pp-1)^{2}\sum_{k=1}^{\lfloor\frac{r-3}{2}\rfloor}(p^{k}-1)= 1-\pp+\pp^{2}-\pp^{\lfloor\frac{r-1}{2}\rfloor}+\pp^{\lfloor\frac{r+1}{2}\rfloor}.$$
Finally, for $0 \leq j < \lfloor\frac{r-1}{2}\rfloor]$ we have
$$\left(M_{\delta}[\lfloor\frac{r-1}{2}\rfloor]+\sum_{k=1}^{\lfloor r-3/2\rfloor} (p^k-1)M_{\delta}[r-1-k]\right)_{\lfloor\frac{r+1}{2}\rfloor}= $$ 
$$
\lfloor \frac{r+1}{2} \rfloor(\pp^{2}-2\pp+1)+\pp-1+(\pp-1)^{2}\sum_{k=1}^{\lfloor\frac{r-3}{2}\rfloor}(\pp^{k}-1)=1-2\pp+\pp^{2}-\pp^{\lfloor\frac{r-1}{2}\rfloor}+\pp^{\lfloor\frac{r+1}{2}\rfloor}.$$
This completes the proof of the claim. 
\end{proof}

\textbf{Step 2:} Let $r \geq 7$. Let $M_{\delta}^{h}$ be the matrix obtained after performing the \underline{column} transformation to $M^{H}_{\delta}$:
\begin{itemize}
    \item $M^{h}_{\delta}[[1]] = M^{H}_{\delta}[[1]]-M^{H}_{\delta}[[2]]$.
\end{itemize}

\begin{claim} \label{claim:matrix2} The matrix $M^{h}_{\delta}$ is a matrix 
of the form 
$$\begin{psmallmatrix}
0 & 1 & 1 & 1 & 1 & \dots & \dots & \dots & \dots& \dots & \dots& 1 \\
 -\pp & \pp^2-\pp+1 & -\pp & 0 & 0 & \dots & \dots & \dots & \dots& \dots & \dots& 0 \\
 \pp & -\pp & \pp^2+1 & -\pp & 0 & \dots & \dots & \dots &\dots & \dots& \dots& 0 \\
 0 & 0 & -\pp & \pp^2+1 & -\pp & 0 & \dots & \dots &\dots & \dots& \dots& 0 \\
\dots & \dots & \dots & \dots &\dots & 0 & \dots  & \dots & \dots& \dots&\dots& 0 \\
 0 & 0 & \dots & -\pp & \pp^2+1 & -\pp & 0 & \dots & \dots& \dots&\dots& 0 \\
0 & \alpha & \dots & \alpha & \alpha+\pp & \beta+2\pp & \pp \left(\pp^{\lfloor \frac{r-1}{2}\rfloor}-1\right) & 0 & \dots& \dots&\dots& 0 \\
 0 & -\alpha & \dots & \dots& -\alpha & -\alpha-\pp & -\pp^{\lfloor\frac{r+1}{2}\rfloor}+\pp^2+1 & -\pp & 0 & \dots &\dots &0\\
  \dots & \dots & \dots & \dots &\dots & \dots & \dots  & \dots & \dots& \dots& \dots &  0 \\
   0 & \dots  & \dots &\dots &\dots &\dots &\dots &\dots & -\pp & \pp^2+1 & -\pp &0  \\
 0 & \dots & \dots & \dots &\dots &\dots &\dots &\dots &\dots & -\pp & \pp^2+1 & -\pp \\
 \sigma(r-1)_{0}-\sigma(r-1)_{1} & \sigma(r-1)_{1} & \dots & \dots & \dots & \dots & \dots& \dots& \dots& \dots& \dots &\sigma(r-1)_{r-1}\\
\end{psmallmatrix}$$ 
where $\alpha= \pp^{\lfloor \frac{r+1}{2}\rfloor}-\pp^{\lfloor \frac{r-1}{2}\rfloor}$, $\beta= \pp^{\lfloor \frac{r+1}{2}\rfloor}-\pp^{\lfloor \frac{r-1}{2}\rfloor}+\pp^2-2 \pp+1$ and $\sigma(r-1)_{0}-\sigma(r-1)_{1}=\pm 1+ \pp g_1(\pp)$ for some $g_{1} \in \Z[x]$.
More formally, we have
\begin{center}
    $M_{\delta}^{H}[i] = \begin{cases}
    (0, 1, \ldots, 1) & \text{if } i=0; \\
    (-\pp, \pp^{2}-\pp+1, -\pp, 0, \ldots, 0) & \text{if } i=1; \\
    (\pp, -\pp, \pp^{2}+1, -\pp, 0, \ldots, 0) & \text{if } i=2; \\
    (0, \ldots, 0, -\pp_{i-1}, (\pp^{2}+1)_{i}, -\pp_{i+1}, 0, \ldots, 0) & \text{if } 2 < i < \lfloor\frac{r-1}{2}\rfloor; \\
    (0, \alpha, \ldots, \alpha, (\alpha+\pp)_{\lfloor\frac{r-3}{2}\rfloor}, (\beta+2\pp)_{\lfloor\frac{r-1}{2}\rfloor}, \pp \left(\pp^{\lfloor \frac{r-1}{2}\rfloor}-1\right)_{\lfloor\frac{r+1}{2}\rfloor}, 0, \ldots, 0) & \text{if } i=\lfloor\frac{r-1}{2}\rfloor; \\
    (0, -\alpha, \ldots ,-\alpha, (-\alpha-\pp)_{\lfloor\frac{r-1}{2}\rfloor},  (-\pp^{\lfloor\frac{r+1}{2}\rfloor}+\pp^{2}+1)_{\lfloor\frac{r+1}{2}\rfloor}, -\pp_{\lfloor\frac{r+3}{2}\rfloor}, 0, \ldots, 0) & \text{if } i=\lfloor\frac{r+1}{2}\rfloor; \\
    (0, \ldots, 0, -\pp_{i-1}, (\pp^{2}+1)_{i}, -\pp_{i+1}, 0, \ldots, 0) & \text{if } \lfloor\frac{r+1}{2}\rfloor < i < r-1; \\
    (\sigma(r-1)_{0}-\sigma(r-1)_{1}, \sigma(r-1)_{1}, \ldots, \sigma(r-1)_{r-1}) & \text{if } i=r-1.
    \end{cases}$
\end{center}
\end{claim}

\begin{proof}
     The proof follows  from Claim \ref{claim:matrix1} after a direct computation. 
\end{proof} 

\begin{remark} \label{rmk:equaldet}
    Notice that since the determinant of a matrix is invariant under row and column transformations, we have $ \text{det}(M_{\delta})= \text{det}(M^{H}_{\delta})=\text{det}(M^{h}_{\delta})$. Hence, verifying $ \text{det}(M_{\delta}) \neq 0$ is equivalent to proving $\text{det}(M^{h}_{\delta})\neq 0$. 
\end{remark}

\begin{example} Take $r=7$. The matrix $M_{\delta}$ is given by  
\[
M_{\delta}=\begin{pmatrix}
 1 & 1 & 1 & 1 & 1 & 1 & 1 \\
 6 \pp^2-11 \pp+5 & 6 \pp^2-10 \pp+5 & 5 \pp^2-9 \pp+4 & 4 \pp^2-7 \pp+3 & 3 \pp^2-5 \pp+2 & 2 \pp^2-3 \pp+1 & (\pp-1) \pp \\
 5 \pp^2-9 \pp+4 & 5 \pp^2-9 \pp+4 & 5 \pp^2-8 \pp+4 & 4 \pp^2-7 \pp+3 & 3 \pp^2-5 \pp+2 & 2 \pp^2-3 \pp+1 & (\pp-1) \pp \\
 4 \pp^2-7 \pp+3 & 4 \pp^2-7 \pp+3 & 4 \pp^2-7 \pp+3 & 4 \pp^2-6 \pp+3 & 3 \pp^2-5 \pp+2 & 2 \pp^2-3 \pp+1 & (\pp-1) \pp \\
 (\pp-1)^2 & (\pp-1)^2 & (\pp-1)^2 & (\pp-1)^2 & \pp^2-\pp+1 & -\pp & 0 \\
 (\pp-1)^2 & (\pp-1)^2 & (\pp-1)^2 & (\pp-1)^2 & (\pp-1)^2 & \pp^2-\pp+1 & -\pp \\
 \sigma(6)_{1} & \sigma(6)_{2} & \sigma(6)_{3} & \sigma(6)_{4} & \sigma(6)_{5} & \sigma(6)_{6} & \sigma(6)_{7} \\
\end{pmatrix}, 
\]
where 
\begin{align*}
    \sigma(6)= & \Big( 4 \pp^4-12 \pp^3+5 \pp^2-\pp+4, 4 \pp^4-12 \pp^3+4 \pp^2+\pp+3, 4 \pp^4-11 \pp^3+4 \pp^2+2 \pp+1,  \\  & \pp \left(4 \pp^3-10 \pp^2+3 \pp+3\right), \pp \left(3 \pp^3-8 \pp^2+2 \pp+3\right), \pp \left(2 \pp^3-5 \pp^2+\pp+2\right), \pp \left(\pp^3-2 \pp^2+2\right)\Big). 
\end{align*} 
After Step 1 we obtain
\[
M_{\delta}^{H}=\begin{pmatrix}
 1 & 1 & 1 & 1 & 1 & 1 & 1 \\
 (\pp-1)^2 & \pp^2-\pp+1 & -\pp & 0 & 0 & 0 & 0 \\
 0 & -\pp & \pp^2+1 & -\pp & 0 & 0 & 0 \\
 (\pp-1) \pp^3 & (\pp-1) \pp^3 & \pp \left(\pp^3-\pp^2-1\right) & \pp^4-\pp^3+\pp^2+1 & \pp \left(\pp^3-1\right) & 0 & 0 \\
 -(\pp-1) \pp^3 & -(\pp-1) \pp^3 & -(\pp-1) \pp^3 & -\pp^4+\pp^3-\pp & -\pp^4+\pp^2+1 & -\pp & 0 \\
 0 & 0 & 0 & 0 & -\pp & \pp^2+1 & -\pp \\
 \sigma(6)_{1} & \sigma(6)_{2} & \sigma(6)_{3} & \sigma(6)_{4} & \sigma(6)_{5} & \sigma(6)_{6} & \sigma(6)_{7} \\
\end{pmatrix}.
\]
After Step 2 we obtain 
\[
M_{\delta}^{h}=\begin{pmatrix}
 0 & 1 & 1 & 1 & 1 & 1 & 1 \\
-\pp & \pp^2-\pp+1 & -\pp & 0 & 0 & 0 & 0 \\
 \pp & -\pp & \pp^2+1 & -\pp & 0 & 0 & 0 \\
 0 & (\pp-1) \pp^3 & \pp \left(\pp^3-\pp^2-1\right) & \pp^4-\pp^3+\pp^2+1 & \pp \left(\pp^3-1\right) & 0 & 0 \\
 0 & -(\pp-1) \pp^3 & -(\pp-1) \pp^3 & -\pp^4+\pp^3-\pp & -\pp^4+\pp^2+1 & -\pp & 0 \\
 0 & 0 & 0 & 0 & -\pp & \pp^2+1 & -\pp \\
 (\pp-1)^2 & \sigma(6)_{2} & \sigma(6)_{3} & \sigma(6)_{4} & \sigma(6)_{5} & \sigma(6)_{6} & \sigma(6)_{7} \\
\end{pmatrix}.
\]

\end{example}

Notice that, from Claim~\ref{claim:matrix2}, to compute the determinant of the matrix $M_{\delta}^{h}$ we need to take a look at the determinant of the $(r-1) \times (r-1)$  submatrix of $M_{\delta}^{h}$ obtained by deleting the $(r-1)$th row and the 1st column; call this submatrix $M$. In particular, $M$ is a Hessenberg matrix (after moving the first row to the last one), that is, it satisfies $M[i][j]=0$ for all $j>i+1$. We will use the following theorem about determinants of Hessenberg matrices to prove Theorem~\ref{prop:maximalrank}. 

\begin{theorem}\cite{fibonacci} \label{thm:HessMat}
    Let $M=(m_{i, j})_{1 \leq i, j \leq m}$ be an $m \times m$ Hessenberg matrix.  Let $M(k)$ denote the submatrix of~$M$ given by $$M(k):=\begin{psmallmatrix}
        m_{1, 1} & \cdots & 0 \\
        \hdots & \ddots & \hdots \\
        m_{1, k} & \cdots & m_{k, k}
    \end{psmallmatrix}.$$ 
    Set $\det(M(0))=1$.  Then the determinanats $\{\det(M(k))\}_{k\geq 0}$ satisfy $\det(M(0))=1$, $\det(M(1))=m_{1, 1}$  and
   $$ \det(M(n)) = m_{n,n} \det(M(n-1))+ \sum_{i=1}^{n-1} \left( (-1)^{n-i}m_{n,i}\left(\prod_{j=i}^{n-1}m_{j, j+1} \det(M(i-1))\right) \right).$$ 
\end{theorem}

Now we are ready to show that $ \text{det}(M^{h}_{\delta})\neq 0$. 

\begin{proposition} \label{thm:detMdelta} We have $\text{det}(M_{\delta}^{h}) = \pm 1 + \pp f(\pp)$ for some $f(x) \in \Z[x]$. 
\end{proposition}

\begin{proof} For $r\leq 6$, a quick computation shows:
    \begin{itemize}
        \item \textbf{Case 1:} $r=2$. We have $\text{det}(M_{\delta}^{h})=1$.
        \item \textbf{Case 2:} $r=3$. We have $\text{det}(M_{\delta}^{h})=1$. 
        \item \textbf{Case 3:} $r=4$. We have $\text{det}(M_{\delta}^{h})= -1$. 
        \item \textbf{Case 4:} $r=5$. We have $\text{det}(M_{\delta}^{h})=1$. 
        \item \textbf{Case 5:} $r=6$. We have $\text{det}(M_{\delta}^{h})=-1 - \pp^2 + \pp^3 + \pp^4 - \pp^5$. 
    \end{itemize}
Hence, the statement is true for $r<7$. Assume $r \geq 7$. Then $M_{\delta}^{h}$ has the form stated in Claim~\ref{claim:matrix2}.
Hence, if we develop the determinant of $M_{\delta}^{h}$ from the first column, from Claim~\ref{claim:matrix2} we have 
$$ \text{det}(M^{h}_{\delta})= \pp (M^{h}_{\delta})_{2, 0}+ \pp (M^{h}_{\delta})_{3, 0} \pm (\sigma(r-1)_{0}-\sigma(r-1)_{1})(M^{h}_{\delta})_{r-1, 0},$$
where $(M^{h}_{\delta})_{i, j}$ denotes the minor of the matrix $M^{h}_{\delta}$ formed after deleting the $i$-th row and $j$-th column. Hence, to prove the result, it is enough to check that 
$$(\sigma(r-1)_{0}-\sigma(r-1)_{1})(M^{h}_{\delta})_{r-1, 0}=\pm 1+ \pp g(\pp)$$ for some $g(x) \in \Z[x]$. By the expression of $\sigma(r-1)_{k}$ computed in the Appendix --  see Remark~\ref{rmk:sigma(r-1)} -- we have that $(\sigma(r-1)_{0}-\sigma(r-1)_{1})=\pm 1+ \pp g_1(\pp)$ for all $r \geq 2$. 
We will now use Theorem~\ref{thm:HessMat} to show that $(M^{h}_{\delta})_{r-1, 1} = 1 + \pp g_{2}(\pp)$ for some $g_{2} \in \Z[x]$. Let $M \coloneqq M^{h}_{\delta}$ and $M(n)$ be the upper left $n \times n$ matrix block of $M$. We proceed by induction on $n$. We claim that for all $n \in \{1, \ldots, r-1\}$ there exists $f_{n}(x) \in \Z[x]$ such that
$$ \text{det}(M(n)) = 1 + \pp f_{n}(\pp).$$

\textbf{Base case:} $n=1$. We have $M(n) = 1-\pp+\pp^{2}$. Hence, the claim is satisfied.

\textbf{Induction step:} Assume the statement is true for all $M(i)$ with $i < n$. Notice the following two facts of the matrix $M$:
\begin{itemize}
    \item $M[i][i]= 1+\pp h_{i}(\pp)$;
    \item $M[i][i+1]= \pp l_{i}(\pp);$
\end{itemize}
for all $0 \leq i \leq n$, where $h_i(x), l_{i}(x) \in \Z[x]$. Hence, by Theorem~\ref{thm:HessMat} we have 
$$ \text{det}(M(n)) = m_{n,n} \text{det}(M(n-1))+ \sum_{i=1}^{n-1} \left( (-1)^{n-i}m_{n,i}\left(\prod_{j=i}^{n-1}m_{j, j+1} \text{det}(M(i-1))\right) \right) = $$
$$(1+\pp h_{n}(\pp))(1 + \pp f_{n-1}(\pp))+ \sum_{i=1}^{n-1} \left( (-1)^{n-i}m_{n,i}\left(\prod_{j=i}^{n-1}(\pp l_{j}(\pp)) (1 + \pp f_{i-1}(\pp))\right) \right) = $$
$$(1+\pp h_{n}(\pp))(1 + \pp f_{n-1}(\pp))+\pp t_{n}(\pp)= 1+\pp g_{n}(\pp),$$
for some $t_{n}(x) \in \Z[x]$. 
Finally, we conclude that  $\text{det}(M(n)) =1 + \pp f_{n}(\pp)$ for all $n$. 

Taking $g_{2}(x)= f_{r-1}(x)$ we get
$$(\sigma(r-1)_{0}-\sigma(r-1)_{1})(M^{h}_{\delta})_{r-1, 0}=(\pm 1+ \pp g_{1}(x))(1+ \pp g_{2}(x))=\pm 1+ \pp g(x),$$ for $g(x) \in \Z[x]$, which proves the proposition.
\end{proof}

As a corollary, we obtain the following. 

\begin{corollary} \label{cor:detMdelta}The matrix $M_{\delta}^{h}$ satisfies $\text{det}(M_{\delta}^{h}) \neq 0$.
\end{corollary}

\begin{proof}
    By Proposition~\ref{thm:detMdelta} we have $\text{det}(M_{\delta}^{h})=\pm 1 + \pp f(\pp)$ for some $f(x) \in \Z[x]$. If $\text{det}(M_{\delta}^{h})=0$, then $\pm 1=\pp f(\pp)$. But then we would get that $\pp$ divides either~$1$ or $-1$, which is not possible because $\pp = q^{\text{deg}(\p)} \geq 3$. Hence, $\text{det}(M_{\delta}^{h})\neq 0$ as we wanted. 
\end{proof}

Now we can prove Theorem~\ref{prop:maximalrank}. 

\begin{proof}[Proof of Theorem~\ref{prop:maximalrank}]
This follows from Corollary~\ref{cor:detMdelta} together with Remark~\ref{rmk:equaldet}.
\end{proof}

\newpage

\section{Appendix} \label{sec:appendix}

\begin{definition}(\cite[3.2]{kevincusp}) \label{def:Dr-1} The following case distinction defines the cuspidal divisor $D_{r-1}$ for all $r \geq 2$ in Theorem~\ref{thm:cuspgen}.
    \begin{enumerate}
    \item If $r = 2$, let $D_{r-1} := C_1$.
    \item If $r\geq 3$ and $r\equiv 3 \mod 4$, let
    \begin{align*}
    D_{r-1} &:= C_{r-1}-(|\p|^r-|\p|^{r-2})C_1\\&+\sum_{2\leq i\leq \frac{r-1}{2}}(|\p|^{r-1} - |\p|^{r-2} - |\p|^{r-2i+1} + |\p|^{r-2i})C_i\\&-\sum_{\lfloor\frac{r+1}{2}\rfloor\leq i\leq r-2}(|\p|^i - |\p|^{\frac{r-1}{2}} + |\p|^{i-\frac{r-1}{2}} - 1)(C_i-|\p|C_{i+1}).
    \end{align*}
    \item If $r\geq 4$ and $r\equiv 0 \mod 4$, let
    \begin{align*}
    D_{r-1} &:= C_{r-1}-(|\p|^r-|\p|^{r-2})C_1\\&+\sum_{2\leq i\leq \frac{r}{2}-1}(|\p|^{r-1} - |\p|^{r-2} - |\p|^{r-2i+1} + |\p|^{r-2i})C_i\\&+\sum_{\substack{\frac{r}{2}\leq i\leq r-2\\\text{$i$: even}}}(|\p|^{i+1}-2|\p|^i+|\p|^{\frac{r}{2}}-|\p|^{i-\frac{r}{2}+1}+1)(C_i-|\p|C_{i+1})\\&-\sum_{\substack{\frac{r}{2}+1\leq i\leq r-3\\\text{$i$: odd}}}(|\p|^{i+1}-|\p|^{\frac{r}{2}}+|\p|^{i-\frac{r}{2}+1}-1)(C_i-|\p|C_{i+1}).
    \end{align*}
    \item If $r\geq 5$ and $r\equiv 1 \mod 4$, let
    \begin{align*}
    D_{r-1} &:= C_{r-1}-(|\p|^r-|\p|^{r-2})C_1\\&+\sum_{2\leq i\leq \frac{r-1}{2}}(|\p|^{r-1} - |\p|^{r-2} - |\p|^{r-2i+1} + |\p|^{r-2i})C_i\\&-\sum_{\substack{\lfloor\frac{r+1}{2}\rfloor\leq i\leq r-2\\\text{$i$: odd}}}(2|\p|^{i+1} - |\p|^i - |\p|^{\frac{r-1}{2}} + |\p|^{i-\frac{r-1}{2}} - 1)(C_i-|\p|C_{i+1})\\&-\sum_{\substack{\frac{r+3}{2}\leq i\leq r-3\\\text{$i$: even}}}(|\p|^i - |\p|^{\frac{r-1}{2}} + |\p|^{i-\frac{r-1}{2}} - 1)(C_i-|\p|C_{i+1}).
    \end{align*}
    \item If $r\geq 6$ and $r\equiv 2 \mod 4$, let
    \begin{align*}
    D_{r-1} &:= C_{r-1}-(|\p|^r-|\p|^{r-2})C_1\\&+\sum_{2\leq i\leq \frac{r}{2}-1}(|\p|^{r-1} - |\p|^{r-2} - |\p|^{r-2i+1} + |\p|^{r-2i})C_i\\&-\sum_{\frac{r}{2}\leq i\leq r-2}(|\p|^{i+1} - |\p|^{\frac{r}{2}} + |\p|^{i-\frac{r}{2}+1} - 1)(C_i-|\p|C_{i+1}).
    \end{align*}
\end{enumerate}
\end{definition}

\begin{proposition}
    \label{rmk:sigma(r-1)} The following table contains the computations of $\sigma(r-1)_{k}$ for all possible $r \geq 2$ and $0 \leq k \leq r-1$. 
\begin{itemize}
    \item \textbf{Case 1:} $r=2$. Then
    $$\sigma(r-1)_{k} = \begin{cases} \pp-1 & \text{if } k=0; \\   
    \pp & \text{if } k=1. \\  
    \end{cases}$$ 
    
    \item \textbf{Case 2:} $r=3$. Then 
    $$\sigma(r-1)_{k} = \begin{cases} -2 \pp^3+3 \pp^2+\pp-2 & \text{if } k=0; \\   
                                 -2 \pp^3+2 \pp^2+\pp-1 & \text{if } k=1; \\  
                                 -\pp^3+\pp^2+\pp & \text{if } k=2. \\  
                                 \end{cases}$$
    
    \item \textbf{Case 3:} $r=4$. 
     $$\sigma(r-1)_{k}  = \begin{cases} -2                         \pp^3+3 \pp-3 &                          \text{if } k=0; \\   
                         -2 \pp^3+\pp^2+3 \pp-2 & \text{if } k=1; \\  
                        -2 + 2 \pp + \pp^2 - \pp^3 & \text{if } k=2; \\  
                        -\pp \left(\pp^2-2\right) & \text{if } k=3. \\  
    \end{cases}$$

    \item \textbf{Case 4:} $r=5$. 
    $$\sigma(r-1)_{k} = \begin{cases} -3 \pp^3+4 \pp^2+2 \pp-3 & \text{if } k=0; \\   
                    -3 \pp^3+3 \pp^2+2 \pp-2 & \text{if } k=1; \\  
                    -2 \pp^3+3 \pp^2+\pp-2 & \text{if } k=2; \\  
                    -2 \pp^3+\pp^2+\pp \left(\pp^2-2\right) & \text{if } k=3; \\  
                    \pp^2 & \text{if } k=4. \\  
    \end{cases}$$

    \item \textbf{Case 5:} $r=6$. 
    $$\sigma(r-1)_{k} = \begin{cases} 
                        -3 \pp^3+3 \pp^2+4 \pp-4 & \text{if } k=0; \\   
                        -3 \pp^3+2 \pp^2+4 \pp-3 & \text{if } k=1; \\  
                        -2 \pp^3+2 \pp^2+3 \pp-3 & \text{if } k=2; \\  
                        \pp^2+3 \pp-2 \pp-2 & \text{if } k=3; \\  
                        -\pp^3+2 \pp-1 & \text{if } k=4; \\ 
                        \pp & \text{if } k=5. \\  
    \end{cases}$$ 
\newpage
\pagestyle{plain}
\newgeometry{bottom=1cm, left=1cm, right=1cm}
\begin{landscape}

\hspace{10cm}
\\

\begin{center}

    \item \textbf{Case 6:} $r\geq 7$ and $r\equiv 3 \pmod{4}$. 
        $$\sigma(r-1)_{k} = \begin{cases} 
            -\frac{(\pp-1)}{2} \left((5-3 r)+(r-1)\pp+2(r-1)\pp^{2}+(r-3)\pp^{\frac{r-3}{2}}-(r+1)\pp^{\frac{r-1}{2}}\right) & \text{if } k=0; \\   
            
            -\frac{(\pp-1)}{2} \left((7-3 r)+(r+1)\pp+2 (r-1)\pp^{2}+(r-3) \pp^{\frac{r-3}{2}}-(r+1)\pp^{\frac{r-1}{2}}\right) & \text{if } k=1; \\  
            
            -\frac{(\pp-1)}{2} \left((7-3 r)+ (r-1)\pp+2(r-2)\pp^{2}+(r-3)\pp^{\frac{r-3}{2}}-(r+1)\pp^{\frac{r-1}{2}}\right) & \text{if } k=2; \\    
            
            \frac{(\pp-1)}{2} \left((-2 k+3 r-3)-(r-1)\pp+2(k-r)\pp^{2}-(r-3)\pp^{\frac{r-3}{2}}+(r+1) \pp^{\frac{r-1}{2}}\right) & \text{if } 2< k< \frac{r-1}{2}; \\ 
            
            -\frac{(\pp-1)}{2} \left(-2(r-1)+(r-1)\pp+(r+1)\pp^2+(r-3)\pp^{\frac{r-3}{2}}-(r+1) \pp^{\frac{r-1}{2}}\right)& \text{if } k=\frac{r-1}{2}; \\ 
            
            -\frac{(\pp-1)}{2} \left(-2(r-3)+(r-1)\pp+(r-1)\pp^{2}+(r-3) \pp^{\frac{r-3}{2}}-(r-1)\pp^{\frac{r-1}{2}}\right) & \text{if } k=\lfloor\frac{r+1}{2}\rfloor; \\
            
            (\pp-1) \left(-2(k-r+1)+(k-r)\pp+(k-r)\pp^{2}+(k-r+1) \pp^{\frac{r-3}{2}}+(r-k)\pp^{\frac{r-1}{2}}\right) & \text{if } \frac{r+1}{2}< k< r-2; \\ 
            
            -2+4\pp-2 \pp^3+\pp^{\frac{r-3}{2}}-3 \pp^{\frac{r-1}{2}}+2 \pp^{\frac{r+1}{2}} & \text{if } k=r-2; \\ 
            2 \pp-\pp^3-\pp^{\frac{r-1}{2}}+\pp^{\frac{r+1}{2}}  & \text{if } k=r-1. \\  
            
    \end{cases}$$

        \item \textbf{Case 7:} $r\geq 8$ and $r \equiv 0 \pmod{4}$. 
            $$\sigma(r-1)_{k} = \begin{cases} 
                        -\frac{(\pp-1)}{2} \left((4-3r)+(r+2)\pp+2(r-2)\pp^2+(r-2)\pp^{\frac{r}{2}-2}-(r+2) \pp^{\frac{r}{2}-1}\right) & \text{if } k=0; \\  
                        
                        -\frac{(\pp-1)}{2} \left(-3(r-2)+ (r+4)\pp+2 (r-2) \pp^2+(r-2)\pp^{\frac{r}{2}-2}-(r+2) \pp^{\frac{r}{2}-1}\right) & \text{if } k=1; \\ 
                        
                        -\frac{(\pp-1)}{2} \left(-3(r-2)+(r+2)\pp+2(r-3)\pp^2+(r-2)\pp^{\frac{r}{2}-2}-(r+2)\pp^{\frac{r}{2}-1}\right) & \text{if } k=2; \\  
                        
                        \frac{(\pp-1)}{2} \left((-2 k+3 r-2)-(r+2)\pp+2(k-r+1)\pp^2-(r-2) \pp^{\frac{r}{2}-2}+(r+2)\pp^{\frac{r}{2}-1}\right) & \text{if } 2<k<\frac{r}{2}-1; \\
                        
                        -\frac{(\pp-1)}{2} \left(-2r+(r+2)\pp+r\pp^2+(r-2) \pp^{\frac{r}{2}-2}-(r+2) \pp^{\frac{r}{2}-1}\right) & \text{if } k=\frac{r}{2}-1; \\
                        
                        -\frac{(\pp-1)}{2} \left(-2(r-1)+r\pp+(r-2)\pp^2+(r-2)\pp^{\frac{r}{2}-2}-r\pp^{\frac{r}{2}-1}\right) & \text{if } k=\frac{r}{2}; \\  
                        
                        (\pp-1)\left(-2(k-r+1)+\pp(k-r-2)+\pp^2(k-r)+(k-r+1) \pp^{\frac{r}{2}-2}+(r-k)\pp^{\frac{r}{2}-1}\right) & \makecell{\text{if } \frac{r}{2} <k< r-2, \\ k=\frac{r}{2}+j \text{ and } j \text{ even};} \\  
                        
                        (\pp-1)\left((-2 k+2 r-1)+(k-r)\pp+(k-r+1)\pp^2+(k-r+1)\pp^{\frac{r}{2}-2}+(r-k) \pp^{\frac{r}{2}-1}\right) & \makecell{\text{if } \frac{r}{2} <k< r-2, \\ k=\frac{r}{2}+j \text{ and } j \text{ odd};} \\ 
                        
                        -(\pp-1) \left(-3+2\pp+\pp^2+\pp^{\frac{r}{2}-2}-2\pp^{\frac{r}{2}-1}\right) & \text{if } k=r-2; \\
                        
                        3 \pp-\pp^2-\pp^3-\pp^{\frac{r}{2}-1}+\pp^{\frac{r}{2}}& \text{if } k=r-1. \\ 
                        

    \end{cases}$$  

    \newpage
    \hspace{12cm}

        \item \textbf{Case 8:} $r \geq 9$ and $r\equiv 1 \pmod{4}$.
            $$\sigma(r-1)_{k}= \begin{cases} 
                        \frac{(\pp-1)}{2} \left(2(r-2)+(r-1)\pp-3 (r-1)\pp^{2}-(r-3) \pp^{\frac{r-3}{2}}+(r+1)\pp^{\frac{r-1}{2}}\right) & \text{if } k=0; \\  
                        
                        -\frac{(\pp-1)}{2} \left(-2(r-3)-(r-3)\pp+3(r-1)\pp^2+(r-3) \pp^{\frac{r-3}{2}}-(r+1) \pp^{\frac{r-1}{2}}\right) & \text{if } k=1; \\

                        -\frac{(\pp-1)}{2} \left(-2(r-3)-(r-1)\pp +(3 r-5)\pp^2+(r-3)\pp^{\frac{r-3}{2}}-(r+1)\pp^{\frac{r-1}{2}}\right) & \text{if } k=2; \\

                        \frac{(\pp-1)}{2}\left(-2(k-r+1)+(r-1)\pp+(2 k-3 r+1)\pp^2-(r-3) \pp^{\frac{r-3}{2}}+(r+1)\pp^{\frac{r-1}{2}}\right) & \text{if } 2 <k< \frac{r-1}{2}; \\ 
                        
                        -\frac{(\pp-1)}{2} \left(-(r-1)-(r-1)\pp +2r\pp^2 +(r-3) \pp^{\frac{r-3}{2}}-(r+1) \pp^{\frac{r-1}{2}}\right) & \text{if } k=\frac{r-1}{2}; \\  
                        
                        -\frac{(\pp-1)}{2} \left(-(r-5)-(r-5)\pp+2(r-1)\pp^2 +(r-3) \pp^{\frac{r-3}{2}}-(r-1) \pp^{\frac{r-1}{2}}\right) & \text{if } k=\frac{r+1}{2}; \\  

                        (\pp-1) \left((-k+r-2)+(-k+r-2)\pp+2(k-r)\pp^{2} +(k-r+1) \pp^{\frac{r-3}{2}}+(r-k) \pp^{\frac{r-1}{2}}\right) & \makecell{\text{if } \frac{r+1}{2} <k< r-2, \\ k=\frac{r+3}{2}+j \text{ and } j \text{ odd};} \\ 
                        
                        (\pp-1) \left((-k+r-1)+(r-k)\pp+(2 k-2 r+1)\pp^{2}+(k-r+1) \pp^{\frac{r-3}{2}}+(r-k) \pp^{\frac{r-1}{2}}\right) & \makecell{\text{if } \frac{r+1}{2} <k< r-2, \\ k=\frac{r+3}{2}+j \text{ and } j \text{ even};} \\  
                        
                        -(\pp-1) \left(4 \pp^{2}+\pp^{\frac{r-3}{2}}-2 \pp^{\frac{r-1}{2}}\right) & \text{if } k=r-2; \\  
                        
                        2 \pp^2-\pp^3-\pp^{\frac{r-1}{2}}+\pp^{\frac{r+1}{2}} & \text{if } k=r-1. \\  
                        

    \end{cases}$$ 

    \item \textbf{Case 9:} $r\geq 10$ and $r \equiv 2 \pmod{4}$.     
    $$\sigma(r-1)_{k} = \begin{cases} 
                        \frac{(\pp-1)}{2} \left(2 (r-2)+ (r-2)\pp+(4-3 r)\pp^2 -(r-2) \pp^{\frac{r}{2}-2}+(r+2) \pp^{\frac{r}{2}-1}\right) & \text{if } k=0; \\   
                        
                        \frac{(\pp-1)}{2} \left(2 (r-3)+(r-4)\pp-3(r-2)\pp^2 +(4-3 r)\pp^2-(r-2) \pp^{\frac{r}{2}-2}+(r+2) \pp^{\frac{r}{2}-1}\right) & \text{if } k=1; \\  
                        
                        \frac{(\pp-1)}{2}\left(2 (r-3)+(r-2)\pp-(r-2) \pp^{\frac{r}{2}-2}+(r+2) \pp^{\frac{r}{2}-1}\right) & \text{if } k=2; \\ 
                        
                        \frac{(\pp-1)}{2} \left(-2 (k-r+1)+(r-2)\pp+(2 k-3 r+2)\pp^2-(r-2) \pp^{\frac{r}{2}-2}+(r+2) \pp^{\frac{r}{2}-1}\right) & \text{if } 2 <k< \frac{r}{2}-1; \\  
                        
                        \frac{(\pp-1)}{2} \left(r+(r-2)\pp-2r\pp^2-(r-2) \pp^{\frac{r}{2}-2}+(r+2) \pp^{\frac{r}{2}-1}\right) & \text{if } k=\frac{r}{2}-1; \\
                        
                        \frac{(\pp-1)}{2} \left((r-2)+(r-4)\pp-2 (r-1)\pp^2-(r-2) \pp^{\frac{r}{2}-2}+r \pp^{\frac{r}{2}-1} \right) & \text{if } k=\frac{r}{2}; \\  
                        (\pp-1) \left((-k+r-1)+(-k+r-2)\pp+(2 k-2 r+1)\pp^2+(k-r+1) \pp^{\frac{r}{2}-2}+(r-k) \pp^{\frac{r}{2}-1}\right) & \text{if } \frac{r}{2} <k< r-2. \\  
                        
                        -(\pp-1) \left(-1+3 \pp^2\pp^{\frac{r}{2}-2}-2 \pp^{\frac{r}{2}-1}\right) & \text{if } k=r-2. \\  
                        
                        \pp+\pp^2-\pp^3-\pp^{\frac{r}{2}-1}+\pp^{\frac{r}{2}}
                        & \text{if } k=r-1. \\  
                        

    \end{cases}$$  
    \end{center}
\end{landscape}

\restoregeometry

\end{itemize}
\end{proposition}

\bibliographystyle{abbrv}
\bibliography{bib}

\end{document}